\documentclass{edpa}
\usepackage[centertags]{amsmath}
\usepackage{amsfonts,amssymb,latexsym,epsfig,amsthm} 
\usepackage{amscd}
\newcommand{\disp}{\displaystyle}
\newcommand{\RR}{\mathbb{R}}
\newcommand{\PP}{\mathbb{P}}
\newcommand{\Z}{\mathbb{Z}}
\newcommand{\X}{\mathbb{X}}
\newcommand{\Y}{\mathbb{Y}}
\newcommand{\Xl}{\mathbb{X}_\ell}
\newcommand{\Xh}{\mathbb{X}_h}
\newcommand{\N}{\mathbb{N}} 
\newcommand{\C}{\mathbb{C}}
\newcommand{\T}{\mathbb{T}}
\renewcommand{\H}{\mathbb{H}}
\renewcommand{\L}{\mathbb{L}}
\newcommand{\EE}{\mathbb{E}} 
\newcommand{\V}{\mathbf{V}}
\newcommand{\bC}{\mathbf{C}}
\newcommand{\one}{\mathbf{1}} 
\newcommand{\dt}[1]{\frac{d#1}{dt}}
\newcommand{\pdt}[1]{\frac{\partial #1}{\partial t}}
\newcommand{\ip}[2]{{\langle#1,#2\rangle}}
\newcommand{\ipX}[2]{{\langle#1,#2\rangle}_\X}
\newcommand{\Ltwo}{\L^2}
\newcommand{\LL}[1]{\left\lvert#1\right\rvert_{\Ltwo}}
\newcommand{\normX}[2][\relax]{\ifx#1\relax\left\lvert#2\right\rvert_{\X}%
        \else\csname#1l\endcsname|#2\csname#1r\endcsname|_{\X}%
        \fi}
\newcommand{\eqdef}{\stackrel{\mbox{\tiny def}}{=}}
\newcommand{\eqM}{\stackrel{\mbox{\tiny M}}{=}}
\newcommand{\bpf}[1][Proof]{{\vspace{1 ex}\noindent {\sc #1:    }}}
\newcommand{\epf}{{{\hspace{4 ex} $\Box$      \smallskip}}}
\newcommand{\ccdot}{ \ \cdot \ }
\newcommand{\SPAN}{\mbox{span}}
\newcommand{\Pl}{\Pi_\ell}
\newcommand{\Ph}{\Pi_h}
\renewcommand{\P}{\mathbf{P}}
\newcommand{\s}{\mathbf{s}}
\renewcommand{\r}{\mathbf{r}}
\renewcommand{\S}{\mathbf{S}}
\newcommand{\R}{\mathbf{R}} 

\newtheorem{assumption}{Assumption}

\newtheorem{theorem}{Theorem}
\newtheorem{lemma}{Lemma}[section]
\newtheorem{corolary}[lemma]{Corollary}

\author{Jonathan C. Mattingly\\{\tiny Original 10/2003, Corrected Version: 2/2004} } 
\title{On Recent Progress for \\the
  Stochastic Navier Stokes Equations}

\address{School of Mathematics, Institute For Advanced Study,
  Princeton NJ, USA 
        and Department of Mathematics, Duke University,
  Durham NC, USA}

\email{jonm@math.duke.edu}

\abstract{We give an overview of the ideas central to some recent
  developments in the ergodic theory of the stochastically forced
  Navier Stokes equations and other dissipative stochastic partial
  differential equations. Since our desire is to make the core ideas
  clear, we will mostly work with a specific example: the
  stochastically forced Navier Stokes equations. To further clarify
  ideas, we will also examine in detail a toy problem. A few general
  theorems are given. Spatial regularity, ergodicity, exponential
  mixing, coupling for a SPDE, and hypoellipticity are all discussed.}
\http{www.math.duke.edu/~jonm} 

\begin{document}
\maketitle

This article attempts to collect a number of ideas which have proven
useful in the study of stochastically forced dissipative partial
differential equations. The discussion will center around those of
ergodicity but will also touch on the regularity of both solutions and
transition densities.  Since our desire is to make the core ideas
clear, we will mostly work with a specific example: the stochastically
forced Navier Stokes equations. To further clarify ideas, we will also
examine in detail a toy problem. Though we have not tried to give any
great generality, we also present a number of abstract results to help
isolate what assumptions are used in which arguments. Though a few
results are presented in new ways and a number of proofs are
streamlined, the core ideas remain more or less the same as in the
originally cited papers. We do improve sightly the exponential mixing
results given in \cite{b:Mattingly02}; however, the techniques used
are the same. Lastly, we do not claim to be exhaustive.  This is not
meant to be an all encompassing review article. The view point given
here is a personal one; nonetheless, citations are given to good
starting points for related works both by the author and others.

Consider the two-dimensional Navier-Stokes equation with stochastic
forcing: 
\begin{equation}
  \label{eq:PureNS}
  \left\{\begin{array}{l}
      \disp\frac{\partial u}{\partial t}
      +(u \cdot\nabla) u +\nabla P=\nu\Delta u+
      \frac{\partial W(x,t)}{\partial t} \\*[.12in]
      \disp\nabla\cdot u = 0
    \end{array}\right. . 
\end{equation}  
We restrict to the $2\pi$-periodic case with mean flow zero, though
many of the results apply equally to bounded domains with Dirichlet
boundary conditions.  The addition of a stochastic forcing can be
motivated by a number of considerations. Since the Navier-Stokes
equations are dissipative, if there is no external forcing, the system
relaxes to the zero state where the fluid is at rest. Hence, if one is
interested in probing the nonlinear dynamics, some forcing is
necessary.  Stochastic forcing is often proposed, particularly in the
study of turbulent fluid flows, as a way to add a ``generic'' forcing.
Generic is then interpreted in the sense of the typical events in
probability space.

We will take the forcing to be the sum of independent Brownian
motions exciting independent Fourier modes. This is convenient
because one of our long term goals is to understand the
interaction between the different scales and the differences
of the dynamics at different scales. Specifically we set
\begin{equation}
  \label{eq:noise}
  W(x,t)= \sum_{k \in \mathcal{K}}
  \sigma_k \frac{k^\perp}{|k|}e^{i k\cdot x} \beta_k(t) \quad
  \mbox{with}  \quad \sum_{k \in \mathcal{K}}
  \sigma_k^2 < \infty
\end{equation}
where $\mathcal{K} \subset \Z^2$ does not contain the zero wave number
ensuring that the spatial mean stays zero. The
$\beta_k=\frac1{\sqrt{2}}(\beta_k^{(1)} + i \beta_k^{(2)})$ where the
$ \beta_k^{(i)}$ are mean zero, variance one Brownian Motions
independent except for the reality condition $\bar
\beta_k=\beta_{-k}$.  The $\sigma_k \in \C$ are constants used to set
the spatial roughness of the flow. They also satisfy the reality
condition $\bar \sigma_k=\sigma_{-k}$. We make the standing assumption
that $\mathcal{E}_0=\sum |\sigma_k|^2 < \infty$ and define
$\sigma_*^2=\max |\sigma_k|^2$. Similarly if $\mathcal{E}_\alpha=\sum
|\sigma_k|^2|k|^{2\alpha} < \infty$\label{er:E1} then for every $t$, $W(\ccdot,t)$ is
almost surely in the Sobolev space $H^\alpha(\T^2) \times
H^\alpha(\T^2)$. Here $\T^2$ is the two dimensional torus. If the
$|\sigma_k|$ decay exponentially or faster, the forcing field is
analytic is space almost surely.

In the next section, we continue with the setup. In section
\ref{sec:invMeasures}, we briefly discuss invariant measures. In
section \ref{sec:FormForcing}, we discuss how the structure of the
solution changes for different choices of forcing. In particular, we
discuss the spatial smoothness. In sections \ref{sec:infiniteDimHard},
\ref{sec:kolmogorovEquation}, and \ref{sec:Rough}, we highlight some
of the difficulties with ergodic theory in infinite dimensions. In
section \ref{sec:effectiveEllipticity}, we discuss ergodicity of the
stochastically forced Navier Stokes (SNS) equations under various
assumptions, including the ideas of ``effective ellipticity'' and the
reduction to Gibbsian dynamics (dynamics with memory). In section
\ref{sec;ergodicity:general}, we formulate the results in a more
general setting and examine a toy model to highlight the main ideas. 
In sections \ref{sec:contractive} and \ref{sec:local}, we discuss the
contractive nature of the SNS equations and the fluctuations of its
energy and enstrophy. In section \ref{sec;genContraction}, we discuss
the Lyapunov structure and localization in the general setting. In
section \ref{sec:expMixing}, we prove a general exponential mixing
result using a non-Markovian coupling argument.
In section
\ref{sec:otherExamp}, we discuss some other systems where the
discussed methods apply. In section \ref{sec:hypoelliptic}, we give a
number of partial results in the setting where the previously stated
ergodic theorems do not hold.  Lastly in section \ref{sec:open}, we
list a few open questions.

\section{The Setting}
\label{sec:setting}
It is convenient to project \eqref{eq:PureNS} onto the space of
divergence free vector fields thereby removing the pressure, which
is just a Lagrange multiplier enforcing the divergence free
constraint. To this end, $\L^2$ will denote the closure in the $L^2$
topology of divergence free, mean zero, $C^\infty$ vector fields on the
two dimensional torus $\T^2$. Similarly the Sobolev space
$\H^\alpha$ is defined as $\L^2$ except that the closure is taken in
$H^\alpha(\T^2) \times H^\alpha(\T^2)$. Projecting equation
\eqref{eq:PureNS} onto $\L^2$ produces the stochastic evolution equation
\begin{align}
  \label{eq:NS}
  \frac{\partial u(x,t)}{\partial t}+ \nu \Lambda^2 u(x,t) +
  B(u,u) = \frac{\partial W(x,t)}{\partial t} 
\end{align}
where $B(u,v)=P_{div}(u\cdot\nabla)v$, $\Lambda^2u=
-P_{div}\Delta u$ and $P_{div}$ is the projection operator onto
the space of divergence free vector fields. 

To better elucidate some of the structure of \eqref{eq:NS}, it is
useful at times to consider the equation for the vorticity
$\omega(x,t)\eqdef\frac{\partial u_2}{\partial x_1} - \frac{\partial
  u_1}{\partial x_2}$ written in Fourier Space. Notice that in two
dimensions $\omega$ is a scalar quantity. Setting $\omega(x,t)=\sum_k
\omega_k(t)e^{i k\cdot x}$, one obtains the infinite system of
coupled diffusions
\begin{equation}\label{eq:vorticity}
  \dt{\omega_k(t)} = -\nu |k|^2 \omega_k + i \sum_{\ell+j=k} \frac{k^\perp
    \cdot \ell}{|k|^2} \omega_\ell \:\omega_j + i|k|\sigma_k\:\one_{k\in
    \mathcal{K}} \: \dt{\beta_k(t)}  \ .
\end{equation}
Unlike many lattices of interacting diffusions, this system in not
invariant under translations in the lattice index $k \in \Z^2$. In
fact for large $|k|$ the linear term in \eqref{eq:vorticity} dominates 
the other drift term which couples the modes together. This
observation is at the heart of all that follows. It gives rise to the
dissipative nature of the dynamics.

Since the noise is additive in our model problem, it is completely
standard that there exists a stochastic flow which depends
continuously on both the initial data and the noise realization $W$
considered as an element of the probability space $\Omega\eqdef
C((-\infty,\infty); \RR^{2|\mathcal{K}|})$. To complete the picture,
we work on $(\Omega,\PP,\mathcal{F},\mathcal{F}_t)$. Here $\Omega$, as
just defined, is the path space of the Brownian trajectories, $\PP$ is
the Weiner measure on this space, $\mathcal{F}$ is the associated
sigma algebra, and $\mathcal{F}_t$ and $\mathcal{F}_{[s,t]}$ are the
filtrations containing the information of the noise increments up to
time $t$ and between time $s$ and $t$ respectively. We will at times
write $\varphi_{s,t}(W)u_0$ or $\varphi_{s,t}^W u_0$ for $u(x,t,W)$
with $u(s)=u_0$ and $\varphi_t(W)$ for $\varphi_{0,t}(W)$. The
notation $u_{[s,t]}$ means the segment of trajectory on $[s,t]$. We
will write $\EE$ to denote expectation with respect to the probability
measure $\PP$; that is $\EE F(W)= \int_\Omega F(W) \PP(dW)$. At times
we will speak of solutions existing on the time interval
$(-\infty,\infty)$. By this we mean a measure $P$ on $\Omega \times
C(-\infty,\infty;\X)$ so that the following holds for almost every
$(W,u)$: $W$ is distributed as a Weiner process, $u(t)$ is adapted to
the filtration generated by $W(s)$ with $s \leq t$, and the pair
$(W,u)$ solves the integral form of \eqref{eq:NS} over any finite time
interval.

\section{Ergodicity and Invariant Measures}
\label{sec:invMeasures}
When investigating a stochastically forced system, such as the
stochastically forced Navier Stokes equation (SNS), the main interest
is often the behavior and structure of the system once it has
forgotten its initial condition. In other words, we are interested in
the behavior of the system in its statistical steady state. The
statistical steady states of a system are described by its invariant
measure. In our setting, a measure $\mu$ on $\L^2$ is invariant under the
dynamics if for any $t>0$ and Borel set $A \subset \L^2$
\begin{align*}
  \mu\{ u_0 : u_0 \in A \} &= \int_\Omega \mu\{ u_0 : u(t,W;u_0) \in
  A\}\PP(dW)= \EE \mu\{ u_0 : \varphi_t(W)u_0 \in A\} \ .
\end{align*}
A system is uniquely ergodic, or simply ergodic, if there is only one
such invariant measure.  The Birkoff ergodic theorem (cf.
\cite{b:Sinai94}) guarantees that for any bounded
function $f:\L^2 \rightarrow \RR$ 
\begin{align*}
  \frac{1}{T} \int_0^T f\Big(u(t,W;u_0)\Big) dt \underset{T\rightarrow
    \infty}{\longrightarrow} \bar{f}(u_0)= \int_{\L^2}
  f(u)d\mu_{u_0}(u)
\end{align*}
if $u_0$ is a typical point for some invariant measure $\mu_{u_0}$. We
have labeled the invariant measure with the initial point $u_0$ to
emphasize that different points might converge to different $\bar f$,
each the average of $f$ against a different invariant measure.
However, if the system is ergodic then there is only one such
invariant measure and the time average $\bar{f}$ is independent of the
initial condition. Hence, the statistics of almost every trajectory
will converge to a unique common distribution. Implying that the
statistics of the systems asymptotic behavior is insensitive to the
initial condition.

\section{The Form of the Forcing}
\label{sec:FormForcing}
Consider the two classes of forcing distinguished by whether
$|\mathcal{K}|< \infty$ or $|\mathcal{K}|=\infty$.  The first class is
the most natural from both the point of view of turbulence theory and
that of exploring the nonlinear dynamics of the Navier-Stokes
equations. In that case, one wants to force the equations at some
scale, usually at large or intermediate scales, and then observe the
transfer of energy and enstropy up and down scale. Generally, forcing
which excites all the Fourier modes ($\mathcal{K}=\Z^2$) is the first
case studied for a given stochastic partial differential equation
(SPDEs). This was true of the SNS (cf.
\cite{b:FlMa95,b:Fe97,b:DaZa96}).  In these investigations, the
forcing was assumed to be spatially rough; essentially $|\sigma_k|
\sim |k|^{-\alpha}$. This assumption means that the forcing is not
analytic in space. The requirement of rough forcing appears to not
simply be a technical assumption, and the methods from
\cite{b:FlMa95,b:Fe97,b:DaZa96} do not seem to work in other elliptic
cases. It is important to mention that the qualitative behavior of the
system appears to be quite different depending on whether the
magnitude of the modes decays at least exponentially or simply
algebraically.

Consider the following two theorems proven respectively in
\cite{b:MattinglySuidan03Pre} and \cite{b:Mattingly02}.  The first
theorem compares the vorticity equation to the associated linear
stochastic heat equation. This equation is just the Ornstein-Uhlenbeck
process
\begin{equation}
  \label{eq:OUpde}
 \left\{ \begin{aligned}
  &\frac{\partial z}{\partial t}(t,x) = \nu\Delta
  z(t,x) + i\sum_k |k|\sigma_k \pdt{\beta_k(t)}\\
  &z(0,x)=\omega_0(x),  
  \end{aligned}\right.  \ .
\end{equation}
If $z(x,t)=\sum_{k \in \Z} z_k(t)\exp(i k\cdot x)$ then
\eqref{eq:OUpde} becomes
\begin{equation}
  \label{eq:OU}
  \dt{z_k(t)}=-\nu|k|^{2} z_k  + i|k|\sigma_k \dt{\beta_k(t)} \ .
\end{equation}
The following theorem states that at small scales $z$ and $\omega$ are
quite similar, even path wise if the forcing decays algebraically in
the spatial Fourier modes.
\begin{theorem}\label{thm:smallScale}
  Assume that $ c|k|^{-\alpha}< |\sigma_k| < C|k|^{-\alpha}$ for some
  positive constants.  Let $\omega_k'=\frac{\sqrt{2} }{|\sigma_k|}
  \omega_k$ and $z_k'=\frac{\sqrt{2} }{|\sigma_k|}
  z_k$.  For any uniformly continuous, bounded function $F$ on $C([0,1];R^d)$,
  $E|F(\omega_{k_1}',...,\omega_{k_d}')-F(z_{k_1}',...,z_{k_d}')|\rightarrow
  0$ as $k_1,...,k_d \rightarrow \infty$. \cite{b:MattinglySuidan03Pre}
\end{theorem}
Theorem \ref{thm:smallScale} says that when the forcing decays algebraically in the
magnitude of the wave number $k$, then so does the solution. In fact, at
small scales, it is pathwise a perturbation of \eqref{eq:OUpde} in
some sense. Hence, the nonlinearity is really secondary in setting the
infinite dimensional character of the problem.

The second theorem covers the case when the forcing decays at least
exponentially fast and, in particular, covers the case when only a
finite number of modes are forced. Earlier versions of this theorem
were proven in \cite{b:Mattingly98b,b:MattinglySinai98} and all of the
versions build on deterministic versions which date back at least to
\cite{b:FoiasTemam89} and are informed by later works such as
\cite{b:Levermore_Oliver97,b:DoeringTiti95,b:DoGi95,b:OliverTiti00}.
In \cite{b:BricmontKupiainenLefevere00,b:Shirikyan02}, yet different
formulations of Theorem \ref{thm:gevry} are given and proven. The
second reference seems to give the best scaling with viscosity, while the
version below gives explicit, eventually stationary processes which
bound the quantities of interest.

\begin{theorem}\label{thm:gevry}
  If there exist positive constants $\beta$ and $C$ so $|\sigma_k| <
  Ce^{-\beta|k|}$ then for any initial $u(0) \in \Ltwo$ there exist two stochastic processes
  $\tau(t,W)$ and $D(t,W)$, positive for $t>0$, so that
\begin{align*}
  |u_k(t,W)| \leq D(t,W)e^{-\tau(t,W)|k|}  
  \text{ $W$-almost surely for all $t> 0$}
\end{align*}
and such that $\lim_{t\rightarrow \infty}\EE\tau(t) \in [c_1,C_1]$ and
$\lim_{t\rightarrow \infty}\EE D(t) \in [c_2,C_2]$ where $c_i$ and
$C_i$ are positive constants which depend on the structure of the
forcing but not on the initial data $u(0)$. (For the form of the
equations for $\tau$ and $D$ and information about their moments see
\cite{b:Mattingly02a}.)
\end{theorem}
Though no lower bound on $|u_k|$, as $|k|\rightarrow \infty$, has been
proven, there is strong evidence that this is the correct order. Even
when the forcing decays faster than exponential, there is no evidence
that the solution does. It is interesting to note that in all of the
current estimates of the decay rate fluctuate in time. Whether this is
correct is not clear. It is a little surprising that even when only a
few modes are forced that $\tau(t)$ does not converge to a constant as
$t \rightarrow \infty$.

Comparing Theorem \ref{thm:smallScale} and \ref{thm:gevry}, one
sees that there is a strong qualitative difference between the two
cases. In the first, the forcing sets the small scale structure. In the
second, the forcing seems to be dictated by the nonlinear dynamics.

\section{The Difficulty of Infinite Dimensions}
\label{sec:infiniteDimHard}
It is reasonable to ask why the ergodic theory of stochastically
forced PDEs is more complicated than that of finite dimensional SDEs.
A basic problem is that there is no single distinguished topology
associated with most infinite dimensional diffusions. Since all
topologies are not equivalent, if one wants to write the transition
density one must use exactly the right base measure. This means one
must know exactly the natural topology of the problem. This is
underlined by the following simple example.  Consider two SPDEs of the
form \eqref{eq:OUpde} with $|\sigma_k|=|k|^{-\alpha}$ in one case and
$|\sigma_k|=|k|^{-\alpha+ \epsilon}$ in the other.  These two process
induce measures on the phase space $L^2$ which are mutually singular
at any positive moment of time, even if they start from the same
point.

In general, getting the correct topology is a very delicate matter.
There seems to be no good general tool to address this class of
problems. In the setting of Theorem \ref{thm:smallScale}, one strongly
suspects that the measure induced by the SNS at a moment of time $t$
is equivalent to that induced by \eqref{eq:OUpde}. However, even in
this case, equivalence has only been proven when the Laplacian is replaced by
$\Delta^{2+\epsilon}$. For this ``hyperviscous'' problem, the
equivalence is proven in \cite{b:MattinglySuidan03Pre}.

\section{Diffusions, Ellipticity, and Hypoellipticity}
\label{sec:kolmogorovEquation}
Just as an ordinary SDE is associated with a PDE which evolves its
density, one can association with an SPDE a ``diffusion'' on a larger
space which evolves the probability transition density. In some cases
this can be made rigorous (cf.
\cite{b:DaZa92,b:FlandoliGozzi98,b:DaPratoZabczyk02}, ). Formally,
consider the ``diffusion'' on $\RR^{Z^2\times Z^2 }$ associated with the
stochastic process \eqref{eq:vorticity}. Writing $z_k=x_k + i y_k$,
the backward Kolmogorov equation would be
\begin{align}\label{eq:BK}
  \frac{\partial \ }{\partial t}U(\{x_k\},\{ y_k\},t)&=
  \mathcal{L}U(\{x_k\},\{y_k\},t)\\ 
  U(\{x_k\},\{ y_k\},0)&=U_0(\{x_k\},\{y_k\}) \notag
\end{align}
where $U_0:\RR^{Z^2\times Z^2}\rightarrow \RR$ is the initial
condition. By $\{x_k\}$ we mean the collection $\{x_k : k \in
\Z^2\}$. The differential operator $\mathcal{L}$ is
\begin{align*} 
  \mathcal{L}= \sum_k \mbox{Re}( F_k) \frac{\partial \ }{\partial
    x_k}+ \mbox{Im}(F_k) \frac{\partial \ }{\partial y_k} +
  \frac12\one_{k\in\mathcal{K}}|k|^2|\sigma_k|^2\big( \frac{\partial^2
    \ }{\partial x_k^2}+ \frac{\partial^2 \ }{\partial y_k^2} \big)
\end{align*}
where 
\begin{align*}
  F_k = -\nu |k|^2 z_k + i
  \sum_{\ell+j=k} \frac{k^\perp \cdot \ell}{|k|^2} z_\ell \:z_j \ .
\end{align*}
The case when $\mathcal{K}=\Z^2$, corresponds to the elliptic setting.
If $|\mathcal{K}|<\infty$, then the operator $\mathcal{L}$ is
degenerate to leading order in all but a finite number of coordinates.
Even the case $\mathcal{K}\not=\Z^2$ but $|\mathcal{K}|=\infty$, it is
still degenerate. In either of the last cases, the ergodic theorems
stated previously are surprising in the sense that they imply some
sort of ellipticity without requiring the detailed geometric
information needed to verify hypoellipticity. These ideas will be
elaborated upon in section \ref{sec:hypoelliptic}.

\section{Ergodicity with Elliptic, Rough  Forcing}
\label{sec:Rough}

In \cite{b:FlMa95,b:Fe97} ergodicity is proven under the assumption,
translated to our setting, that $c|k|^{-\alpha} < |\sigma_k|<
C|k|^{-\alpha}$ for some positive constants. The proof of ergodicity
relies on the Bismuth-Elworthy-Li formula and seems to fundamentally
require an elliptic diffusion with algebraically decaying spectrum. In
light of Theorem \ref{thm:smallScale}, it is tempting to characterize
the system in this regime as a perturbation of the linear process
since the linear process sets the small scale structure. Eckmann and
Hairer \cite{b:EckmannHairer01b} showed that finite dimensional
Malliavin calculus could be combined with the type of analysis used in
\cite{b:Cerrai99,b:FlMa95,b:Fe97} to show that a stochastically forced
SPDE was ergodic even if a finite number of the directions with
possibly positive Lyapunov exponents were not forced. They required a
bracket condition in the spirit of H\"oromander's ``sum of squares
theorem'' (see section \ref{sec:hypoelliptic}). Unfortunately they
still required rough (algebraically decaying) forcing.

\section{Ergodicity under an Effective Ellipticity Assumption}
\label{sec:effectiveEllipticity}
We now turn to a number of results which allow one to prove ergodicity
despite the fact that $\mathcal{K}\not=\Z^2$. In particular, no lower
bound will be placed on the decay rate of the $|\sigma_k|$; even
$|\mathcal{K}| < \infty$ will be allowed if other assumptions are
satisfied. Recalling that $\mathcal{E}_0=\sum_k|\sigma_k|^2$, we have the
following theorem.
\begin{theorem}
  \label{thm:mainErgodic} 
  There exists a fixed constant $\mathcal{C}$ depending only on the
  domain so that the following hold:
  \begin{itemize}
  \item If $\mathcal{C}\frac{\mathcal{E}_0}{\nu^3} < 1$ then
    \eqref{eq:NS} has a unique  $\L^2$-valued invariant probability
    measure regardless of the structure of the forcing.
    \cite{b:Mattingly98b,b:Mattingly98}
  \item  If $|\sigma_k|>0$ for all $k$ with $|k|^2
    \in(0,\mathcal{C}\frac{\mathcal{E}_0}{\nu^3})$, then  \eqref{eq:NS}  has
    a unique $\L^2$-valued invariant probability  measure.
    \cite{b:EMattinglySinai00,b:BricmontKupiainenLefevere01}
  \end{itemize}
\end{theorem}
By a $\L^2$-valued probability measure, we mean a measure $\mu$ such
that $\mu(\L^2)=1$.  The existence was given in
\cite{b:VishikFursikov88,b:Fl94} in the case of the SNS and in a more
general setting in \cite{b:ChowKhasminskii98}. Both results of Theorem
\ref{thm:mainErgodic} stem from the following fact first proven in the
stochastic setting in \cite{b:Mattingly98b} but closely related to
ideas in
\cite{b:FoiasProdi67,b:Temam95,b:ConstantinFoiasNicolaenkoTemam89,
  b:FoiasSellTemam88}.  Contemporaneously to
\cite{b:EMattinglySinai00} similar techniques were used in
\cite{b:KuksinShirikyan00}, to prove a similar theorem for impulsive
or ``kicked'' forcing. Though these initial results applied only for
bounded forced, those authors later extended them to cover unbounded forcing.
They also proved a convergence theorem of the kicked case to the white
in time case.  For the remainder of the discussion of the SNS, we fix
a positive $N_*$.  Let $\Pl$ be the orthogonal projection onto the
space spanned by the wave numbers $k$ with $|k|<N$ and let $\Ph$ be
the complimentary orthogonal projection.  We consider the ``high
mode'' equation on $\Ph\L^2$ given by
\begin{equation}
  \label{eq:highMode}
  \pdt{h(x,t)}+ \nu \Lambda^2 h(x,t) +
  \Ph B(h+\ell,h+\ell)  = \frac{\partial \eta(x,t)}{\partial t}
\end{equation}
where $\ell$ is a given ``low mode'' trajectory in $\Pl\L^2$ and
$\eta=\Ph W(x,t)$. We will denote by
$\Phi_{s,t}^\eta(\ell_{[s,t]};h_0)$ the solution to
\eqref{eq:highMode} at time $t$ with initial condition $h_0$ at time
$s$ and the given external forcings $\ell$ and $\eta$ over the time
interval $[s,t]$. A more quantitative version of the following result
is given in Lemma \ref{l:abs_contraction}.
\begin{theorem}[Foias and Prodi'67, Mattingly'98]\label{thm:determiningModes}
  Let $\mathcal{C}$ be the same constant as in Theorem
  \ref{thm:mainErgodic}. Assuming that $N_*^2 \geq
  \mathcal{C}\frac{\mathcal{E}_0}{\nu^3}$, there exists a positive
  constant $\gamma$ so the following two
  statements hold.
  \begin{itemize}
  \item Let $u(x,t,W)$ be a solution to \eqref{eq:NS} on the time
    interval $[0,\infty)$. Define $\ell(t)=\Pl u(t)$ and
    $\eta(t)=\Ph W(t)$. For almost every $W$, there exists a positive
    constant $T=T(W,u(0))$ so that for all $t \geq T$ and $h_0 \in \Ph\L^2$
    \begin{align*}
      \LL{\Phi_{0,t}^\eta(\ell_{[0,t]};h_0)-\Ph u(t,W)} \leq
      \LL{h_0-\Ph u(0)}e^{-\gamma t} \ .
    \end{align*}
    
    In particular, if $\tilde u(x,t, W)$ is another solution on
    $[0,\infty)$ and $\Omega_0 \subset \Omega \times \Omega$ such that
    for all $(W,\tilde W) \in \Omega_0$ and $t \in [0,\infty)$ one has
    $\Pl \tilde u(t,W)= \Pl u(t,\tilde W)$ and $\Ph W(t)-\Ph W(0)= \Ph
    \tilde W(t)- \Ph \tilde W(0)$ then $\tilde u(t,W)= u(t,\tilde W)$
    for all $t \in[0,\infty)$ and almost every $(W,\tilde W) \in
    \Omega_0$.

  \item Let $u(x,t,W)$ be a stationary solution to \eqref{eq:NS} on
    the time interval $(-\infty,\infty)$. Define $\ell(t,W)=\Pl
    u(t,W)$ and $\eta(t)=\Ph W(t)$. Then with probability one, there
    exists a positive constant $C$ depending only the
    solution $u$ so that for $t \leq 0$
    \begin{align*}
      \LL{\Phi_{t,0}^\eta(\ell_{[t,0]};h_0)-\Ph u(0)} \leq
      C(\LL{h_0} + 1) e^{-\gamma |t|} .
    \end{align*}
    
    In particular, if $\tilde u(x,t, W)$ is another stationary solution
    on $(-\infty,\infty)$, $\Omega_0 \subset \Omega \times \Omega$,
    and $T$ a fixed time, such that for any $(W,\tilde W) \in
    \Omega_0$ and $s \in(-\infty,T]$, $\Pl \tilde u(s,\tilde W)= \Pl
    u(s,W)$ and $\Ph W(s)- \Ph W(0) = \Ph \tilde W(s)- \Ph \tilde
    W(0)$ then $u(s,W)= \tilde u(s,\tilde W)$ for all $s
    \in(-\infty,T)$ and almost every $(W,\tilde W) \in \Omega_0$.
    
    In other words, the history of the modes with wave number $|k|$
    less than $N_*$ combined with the history of the forcing
    increments on the remaining degrees of freedom is sufficient to
    determine the solution uniquely with probability one.
  \end{itemize}
\end{theorem}

The first statement in Theorem \ref{thm:mainErgodic} is really a
consequence of the contractive properties used to prove Theorem
\ref{thm:determiningModes}. It is the special case when the set of
determining low modes is empty; hence, knowledge of the infinite past
of the random forces is sufficient to reconstruct the state of the
whole system.  In general, as shown in Theorem
\ref{thm:determiningModes}, one needs some finite number of
determining modes and knowledge of the random forcing applied to the
missing modes to reconstruct the missing modes.  

We now give a more general result which implies the first part of
Theorem \ref{thm:mainErgodic} by showing that to each realization of
noise there corresponds a unique, stationary solution if the viscosity
is large enough relative to the forcing. Another way of saying this is
that the system's random attractor, whose existence was proven at any
viscosity by Flandoli \cite{b:Fl94}, consists of a trivial diffusing
point.  Schmalfuss proved a similar statement using a random fixed
point argument in the case of multiplicative noise and large viscosity
\cite{b:Schmalfuss97}. In that case, the attracting random solution is
a random fixed point which does not fluctuate in time.

\begin{theorem}\label{thm:largeNu}
  If $\mathcal{C}\frac{\mathcal{E}_0}{\nu^3} < 1$, then there exists a
  unique stationary random solution $u^*(t,W)$ defined for $t \in
  (-\infty,\infty)$ and almost all $W\in \Omega$.  In addition, it
  attracts all other solutions exponentially quickly.
  \cite{b:Mattingly98b,b:Mattingly98}
\end{theorem}

One of the interesting interpretations of Theorem
\ref{thm:determiningModes} in the case of arbitrary viscosity is that
on the set of stationary solutions one can define a functional
$\Phi:C(-\infty,0;\Pl\L^2) \rightarrow \Ph\L^2$ which reconstructs the
high modes from the low modes. In particular if $u$ is a stationary
solution and $\eta=\Ph W$ then define
\begin{align*}
  \Phi^\eta(\Pl u_{(-\infty,0]})\eqdef\lim_{t\rightarrow -\infty}
  \Phi_{t,0}^\eta(\Pl u_{(t,0]};h_0)
\end{align*} 
for some arbitrary fixed $h_0$. Theorem \ref{thm:determiningModes}
guarantees that the limit exists, that it is independent of the choice
of $h_0$, and that $\Ph u(0,W) = \Phi^\eta(\Pl u_{(-\infty,0]})$. With
this result, we can close the low mode equations at the price of
introducing memory. One obtains
\begin{equation}
  \label{eq:lowModeMemory}
  \dt{\ell(x,t)}+ \nu \Lambda^2 \ell(x,t) +
  \Pl B(\ell+\Phi^{\theta_t\eta}(\theta_t\ell),\ell+
  \Phi^{\theta_t\eta}(\theta_t\ell))=  
  \frac{\partial \xi(x,t)}{\partial t}  
\end{equation}
where $\theta_t$ is the shift defined on $\ell$ by
$(\theta_t\ell)(s)=\ell(s+t)$ and $\eta$ by
$(\theta_t\eta)(s)=\eta(t+s)-\eta(t)$.  This representation is closely
related and inspired by the inertial form representation from inertial
manifolds theory (cf.  \cite{b:ConstantinFoiasNicolaenkoTemam89,
  b:EdFoNiTe94}) and the ideas of symbolic dynamics.  From the
representation in \eqref{eq:lowModeMemory}, it is clear why it might
be reasonable to call systems satisfying the assumptions of the
second part of Theorem \ref{thm:mainErgodic} ``effectively elliptic''
diffusions.  Under that assumption, the system reduces to an equation
of the form \eqref{eq:lowModeMemory}. This It\^o process with memory
is elliptic in the sense that the noise directly agitates all of the
coordinates.  In contrast to the hypoelliptic systems considered in
section \ref{sec:hypoelliptic}, no detailed knowledge of the tangent
space structure is needed. Once the assumption about all of the
possibly unstable directions being forced is satisfied, only some soft
general estimates are needed.

When viewed in the context of Section \ref{sec:kolmogorovEquation},
Theorem \ref{thm:mainErgodic} might seem surprising. The theorem
allows the associated diffusion to be degenerate in an infinite number
of directions; yet the system has nice ergodic properties. Yet in other
ways, Theorem \ref{thm:mainErgodic} is expected. It simply says that if
all of the unstable directions are forced directly, the system is
ergodic. Since the long time dynamics are governed by the behavior on the
``unstable manifold'' (if one was known to exist), forcing those
directions destroys all possible obstruction to mixing in the phase
space. Since these systems are non-autonomous, when we say that a
collection of directions are stable, we really mean that all of the
associated Lyapunov exponents associated with these degrees of freedom
are negative. 

\section{Ergodicity: General Constructions}
\label{sec;ergodicity:general}
We now lay out a more general framework to make some ideas clear
without being encumbered by specifics. In the next
section, we also give a simple toy model and some illustrative
examples which hopefully will make the ideas concrete.

Let $(\X,\normX{\ccdot},\ipX{\ccdot}{\ccdot})$ be a complete separable
Hilbert space with a basis $\{e_k\}$, $k=1,\dots$ . Consider the
stochastic evolution equation
\begin{equation}
  \label{eq:abs_equation}
  \dt{u} = G(u) + \dt{W} \ .
\end{equation}
taking values in $\X$. Let $\mathcal{D}(G) \subset \X$ be the domain
of $G$.\label{er:G1} For concreteness, we take $W(x,t)=\sum \sigma_k e_k b_k(t)$
where $\sigma_k$ are constants which fix the structure of the forcing
and the $b_k(t)$ are standard variance one Brownian motions. More
general forcings built over a cylindrical Wiener space are possible
with further assumptions, but this will be sufficient for our needs.

We assume that the $\sigma_k$ are chosen so that
\eqref{eq:abs_equation} has a globally defined stochastic flow
$\varphi^W_{s,t}u_0=u(t,W)$ where $u(s)=u_0$.  It is standard to
associate with this flow a random dynamical system defined by the skew
flow $\Theta_t(u,W)=(\varphi^W_{0,t}u_0,\theta_t W)$ (cf.
\cite{b:Arnold98,b:Kifer86}). Here $\theta_t$ is the shift operator.
On noise paths the shift is defined by $(\theta_t W)(s)=W(t+s)-W(t)$.
We also define the shift of a trajectory by $(\theta_t u)(s)=u(t+s)$.
The difference in definition is due to the fact that in the first case
we are really shifting the noise increments and not the path itself.

Fix a positive integer $N_*$, and define the splitting of the space
$\X=\Xl \times \Xh$, by $\Xl=\SPAN\{e_k: k < N_*\}$ and
$\Xh=\SPAN\{e_k: k \geq N_*\}$. Let $\Pl$ and $\Ph$ be the orthogonal
projectors onto $\Xl$ and $\Xh$ respectively. We will write
$u=(\ell,h)=(\Pl u,\Ph u) \in \Xl \times \Xh$ and $\eta=\Ph W$ and
$\xi=\Pl W$. Notice that the probability measure $\PP$ decomposes into
$\PP_\eta \times \PP_\xi$. As before, we will denote segments of
trajectories by an interval of time as a subscript. Hence,
$\ell_{[s,t]}$ is a trajectory in $\Xl$ between time $s$ and $t$. We
use $\Pi_{[s,t]}$ to denote the projection of a path or set of paths
onto the time interval $[s,t]$.

One can always split the system into two coupled equations on $\Xl
\times \Xh$, 
\begin{align}
  \label{eq:abs_h}
  \dt{h} &= \Ph G(\ell+h) + \dt{\eta}\\
  \dt{\ell} &= \Pl G(\ell+h) + \dt{\xi} \ . 
 \label{eq:abs_l} 
\end{align}
As in section \ref{sec:toy}, given this splitting, one can usually define a
map $h(t)=\Phi_{s,t}^\eta(\ell_{[s,t]};h_0)$ which solves
\eqref{eq:abs_h} given an initial condition $h_0$, noise path $\eta$,
and $\ell_{[s,t]}$ viewed as an external input.\label{er:Phi1} Then for each $t_0$,
$h_0$, and $\eta$, we can define
\begin{align}
  \label{eq:abs_reduced}
   \dt{\ell} &= \Pl G(\ell+\Phi_{t_0,t}^\eta(\ell_{[t_0,t]};h_0)) +
   \dt{\xi}\\
   \ell(t_0)&=\ell_0  \ . \notag
\end{align}
Equation \eqref{eq:abs_reduced} is no longer a standard diffusion as
we have introduced memory through the function $\Phi_{t_0,t}^\eta$. It
is critical to notice that $\ell(t)$ remains an adapted It\^o process
and hence the power of stochastic calculus can be brought to bear.

For the representation in \eqref{eq:abs_reduced} to be useful in the
study of the ergodic theory of \eqref{eq:abs_equation}, the reduced
dynamics \eqref{eq:abs_reduced} must ``forget'' the choice of $h_0$.
One way to investigate this is to study the system as $t_0 \rightarrow
-\infty$.  If the functional $\Phi$ becomes independent of $h_0$, then
we have a closed dynamics on $C(-\infty,0;\Xl)$ over the probability
space $\Omega$. The resulting stochastic process could have infinite
memory.  Since it is defined by a compatible family of Gibbs measures,
in \cite{b:EMattinglySinai00} it was dubbed ``Gibbsian dynamics'' to
be contrasted with Markovian dynamics. The ergodic theory of systems
with this type of memory was explored in its own right in
\cite{b:Bakhtin02,b:BakhtinMattingly03Pre}.

Alternatively, one could study the measures induced on the infinite
future for different choices of $h_0$ and show that they induce the
same asymptotic dynamics in some sense. This was the point of view
taken in \cite{b:Mattingly02}. 

The two approaches are more or less equivalent and each has its own
difficulties. One difficulty of the memory/Gibbsian Dynamics approach
is that sometimes the limit, $\lim_{t_0\rightarrow -\infty}
\Phi_{t_0,t}^\eta(\ell_{[t_0,t]}; h_0)$, only exists on a restricted
set of paths. In any situation where the approach works, one can
always take $\ell_{(-\infty,t]}$, which are typical realizations of a
stationary solution obtained by suspending any invariant measure over
path space. But such a characterization is not constructive and at
times is difficult to work with.

At the most basic level, the success of the approach developed in
\cite{b:EMattinglySinai00} (or
\cite{b:KuksinShirikyan00,b:BricmontKupiainenLefevere01} for that
mater) hinges on treating the $\ell$ and $h$ variables in
fundamentally differently way. Since the $\ell$ variable is finite
dimensional in all the situations we consider, all the difficulties of
probabilistic calculations in an infinite dimensional setting,
mentioned in section \ref{sec:infiniteDimHard}, are not an issue. In
particular, the time $t$ transition densities projected onto $\Xl$
will have densities relative to Lebesgue measure on $\Xl$ if all of
the directions in $\Xl$ are forced. The analysis of the $h$ variable
is dynamic in nature. The analysis is done noise realization by noise
realization.  In contrast the analysis of the $\ell$ variable is
probabilistic in nature. Arguments are made at the level of transition
densities. If the system is strongly contractive, then the structure
of the forcing is irrelevant.  This was the fundamental fact used in
\cite{b:Mattingly98} to prove ergodicity by showing the existence of a
distinguished globally attracting solution. Another way to say this is
that the random attractor is trivial, consisting of a single point at
each moment of time. Given our splitting, a similar structure remains
in the $h$ variable. As we will see, such contraction, $\eta$-fiber by
$\eta$-fiber, is much less sensitive to the topology than are questions
like the absolute continuity of measures. The basic idea is to change
the measure on the $\ell$ variables in such a way that the remaining
degrees of freedom are contractive. The analyses  in
\cite{b:Mattingly02} and \cite{b:EMattinglySinai00} accomplish this by
making the $\ell$'s agree after some finite time.  In
\cite{b:Hairer02}, the measure is changed to bring the $\ell$ (and
$h$) together asymptotically at infinity but never at a finite time.
In all cases, care must be taken so that the changes in the measure to not
accumulate to the extent that the limiting measures become singular.

To execute this program, we need to analyze the dynamics on the path
space of $\Xh$ and understand the structure of the measures induced on
the path space of $\Xl$. To this end, we make a few definitions. For 
all $t>s \geq 0$ define 
\begin{align}
  \label{eq:defQ}
  Q_t(\ell_0,h_0,A)&=\PP\big( \ell(t) \in A \big| \ell(0)=\ell_0,
  h(0)=h_0\big) \\Q_{[s,t)}(\ell_0,h_0,B)&=\PP\big( \ell_{[s,t)} \in B
  \big| \ell(0)=\ell_0, h(0)=h_0\big) \notag
\end{align}
for Borel sets $A\subset \Xl$ and $B \subset C([s,t); \Xl)\cong
C([0,t-s); \Xl)$.  Notice we have associated $C([s,t),\Xl)$ with
$C([0,t-s),\Xl)$ and will view $u_{[s,t)}$ as an element of
$C([0,t-s),\Xl)$.

Similarly for any realization of $\eta$, let
$\mathcal{F}^\eta_{[s,t]}$ be the $\sigma$-algebra generated by the
increments of $\eta$ between $[s,t]$. We define
$Q_t^\eta(\ell_0,h_0,A)=\PP( \ell(t) \in A | \ell(0)=\ell_0, h(0)=h_0,
\mathcal{F}^\eta_{[0,t]})$ and $Q_{[s,t)}^\eta(\ell_0,h_0,B)=\PP(
\ell_{[s,t)} \in B | \ell(0)=\ell_0, h(0)=h_0,
\mathcal{F}^\eta_{[s,t]})$. These are analogous to the previous
measures except that we have conditioned on the realization of $\eta$
over the time interval in question. Hence, for $A \subset \Xl$.
\begin{align*}
  Q_t(\ell_0,h_0,A)=\int Q_t^\eta(\ell_0,h_0,A)\PP(d\eta) = \EE
  Q_t^\eta(\ell_0,h_0,A) \ .
\end{align*}

\subsection{A Toy Problem}
\label{sec:toy}

We now describe a simple toy problem which contains the main ideas
needed to prove the results of the previous section. We will use the
same notation to make the connections explicit.

Consider the following two dimensional stochastic differential
equation
\begin{align}
  \label{eq:toy_h}
  \dt{h(t)} &= -\nu_1 h + F_1(\ell,h) + \sigma_1\dt{\eta}\\
  \dt{\ell(t)} &= -\nu_2 \ell + F_2(\ell,h) + \sigma_2\dt{\xi} .
  \label{eq:toy_l}
\end{align}
Here $\nu_i> 0$, $\sigma_i\geq 0$, $\eta$ and $\xi$ are standard one
dimensional Brownian Motions on the probability space
$\Omega=C((-\infty,\infty);\RR^2)$. Hence, in the notation of the
previous section $\X=\RR^2$, $\Xl=\RR$, and $\Xh=\RR$. We assume the
following estimates hold $|F_1| + |F_2| < K$ and $|F_i(\ell,h)-
F_i(\ell,\tilde h)| \leq L_i |h-\tilde h|$.  For the moment, we allow
either or both of the $\sigma_i$ to be zero.  Eventually, we will
require only that $\sigma_2>0$ allowing $\sigma_1$ to be zero if
desired. Since the $F_i$ are uniformly bounded, it is easy to see that
$\limsup_{t\rightarrow \infty} \EE[ h^2(t) + \ell^2(t)]$ is uniformly
bounded over all initial conditions.  From this, one can deduce the
existence of an invariant measure using standard tightness arguments.
The stochastic flow $\varphi_t^{(\xi,\eta)}(\ell_0,h_0)$ and the
functional $\Phi^\eta_{s,t}(\ell_{[s,t]};h_0)$ are defined as in the
previous section.

Subtracting two copies of \eqref{eq:toy_h} with the same $\eta$ and
$\ell_{[s,t]}$ but different initial conditions produces the estimate
\begin{align}\label{eq:toyContractive}
  |\Phi^\eta_{s,t}(\ell_{[s,t]};h_0) - \Phi^\eta_{s,t}(\ell_{[s,t]};\tilde
   h_0)| \leq |h_0 - \tilde h_0| e^{-(\nu_1 - L_1)(t-s)} \ .
\end{align}
Using this estimate immediately produces the following result, which
is the analog of Theorem \ref{thm:determiningModes}.
\begin{lemma} \label{l:toy_contraction}Assume $\nu_1 > L_1$. Given
  $\ell \in C((-\infty ,0];\RR)$ and $h_0 \in \RR$, the limit
  \begin{equation*}
    \Phi^\eta(\ell_{(-\infty,0]})\eqdef \lim_{s \rightarrow - \infty} 
    \Phi^\eta_{s,0}(\ell_{[s,0]};h_0) 
  \end{equation*}
  is well defined almost surely and independent of $h_0$. Similarly,
  fixing a time interval $[s,t]$ and initial conditions
  $(\ell_i(s),h_i(s))$. Let $\Omega_0 \subset \Omega \times \Omega$
  such that for all $(\xi_1,\eta_1,\xi_2,\eta_2) \in \Omega_0$
  $\eta_1(r)-\eta_1(s)= \eta_2(r)-\eta_2(s)$ and $\ell_1(r)=\ell_2(r)$
  when $r \in [s,t]$ where
  $(\ell_i(r),h_i(r))=\varphi_t^{(\xi_i,\eta_i)}(\ell_i(s),h_i(s))$.
  Then for all $(\xi_1,\eta_1,\xi_2,\eta_2) \in \Omega_0$
  \begin{align*}
    |h_1(t) - h_2(t)| \leq  |h_1(s) - h_2(s)|e^{-(\nu_1 -L_1)(t-s)} \ . 
  \end{align*}
\end{lemma}
Recalling that the shift $\theta_t$ on trajectories acts by
$(\theta_t\ell)(s)= \ell(t+s)$ and on noise paths by
$(\theta_t\eta)(s)=\eta(t+s)-\eta(t)$, then by Lemma
\ref{l:toy_contraction} we can reduce the system to the following
system with memory
\begin{align}\label{eq:toy_mem}
  \dt{\ell} &= -\nu_2 \ell + F_2(\ell,
  \Phi^{\theta_t\eta}(\theta_t\ell)) + \sigma_2\dt{\xi}\\
  \ell(0) &=\ell_0 \notag
\end{align}
where now $\ell(t)$ is seen as an element of $C((-\infty,t];\RR)$.
Similarly our initial condition $\ell_0$ is an element of
$C((-\infty,0];\RR)$.

We now turn to another auxiliary result which, along with the
contraction embodied Lemma \ref{l:toy_contraction}, is the linchpin on
which ergodicity hangs.  Recalling the definitions from
\eqref{eq:defQ}, we have
\begin{lemma}\label{l:toy_marginals} Assume $\sigma_2>0$ and $\nu_1
  >L_1$. For all $\ell_0,h_0 \in \RR$, the measure
  $Q_t(\ell_0,h_0,\ccdot)$ is equivalent to Lebesgue measure.  For all
  $\ell_0,h_0,\tilde h_0, \in \RR$, the measure
  $Q_{[0,\infty)}(\ell_0,h_0,\ccdot)$ is equivalent to
  $Q_{[0,\infty)}(\ell_0,\tilde h_0,\ccdot)$. For any realization of
  $\eta$, the exact same conclusions hold with $Q_t$ replace by
  $Q_t^\eta$ and $Q_{[0,\infty)}$ replaced by $Q_{[0,\infty)}^\eta$.
\end{lemma}

In the next section, we will use the Lemma \ref{l:toy_marginals} to
prove the ergodicity of the toy problem (equations \eqref{eq:toy_h}
and \eqref{eq:toy_l}). Of course, if $\sigma_1,\sigma_2 >0$ then the
system is uniformly elliptic and the fact that there is a unique
invariant measure follows from standard elliptic theory. Even when
$\sigma_1=0$, one might well use hypoelliptic diffusion theory to
prove ergodicity. What we present here is a different possible route,
where the detailed knowledge of the tangent space structure used in
hypoelliptic arguments is replaced with assumptions about the system's
Lyapunov exponents. The advantage of this route being that the contractive
properties are less sensitive to the choice of topology than the
measure theoretic properties of the system needed for the more standard
approaches to ergodicity.

The fact that $Q_{[0,\infty)}(\ell_0,h_0,\ccdot)$ and
$Q_{[0,\infty)}(\ell_0,\tilde h_0,\ccdot)$ are equivalent measures on
the infinite time $[0,\infty)$ interval is critical. Absolute
continuity on finite time intervals would not be sufficient. As an
illustrative example consider the measures induced on path space by a
standard Brownian motion $B(t)$ and the SDEs
\begin{xalignat*}{2}
  \dt{X(t)}&=-X(t) + \dt{B(t)}  & \dt{Y(t)}&=-Y(t) +\frac1t + \dt{B(t)}\\
  X(0)&=x_0 & Y(0)&=x_0 \ .
\end{xalignat*}
All three processes induce measures which are pairwise equivalent on
any finite segment of path space. However, only the processes $X$ and
$Y$ are equivalent on the infinite futures because their difference,
$\frac1t$, is square integrable on an infinite time interval. See the
proof of the second part Lemma \ref{l:toy_marginals} for the needed
argument. In particular, we see that $X(t)$ and $Y(t)$ have the same
asymptotic behavior at the level of the path space marginals, while
$W(t)$ has a different one. Notice that we do not mean that
$|X(t)-Y(t)| \rightarrow 0$ as $t \rightarrow \infty$.

Intuitively it is clear why Lemma  \ref{l:toy_marginals}  when
combined with the contractive estimate from \eqref{eq:toyContractive},
implies that there is only one invariant measure. From Lemma
\ref{l:toy_marginals}, we see that any two invariant measures will
induce equivalent measure in $C(0,\infty;\Xl)$. Hence they will charge
trajectories with the same projection onto $\Xl$. This is already
enough to ensure that the distribution on $\Xl$ is unique. However
because of  \eqref{eq:toyContractive} if the two paths share the same
projection on to $\Xl$ for all time the remaining degrees of freedom
will also converge. Hence the time averages along some typical paths of the
two measure will be the same. This implies the measure are the
same. In Theorem \ref{thm:ergodic}, we make this argument precise.

\bpf[Proof of the first part of Lemma \ref{l:toy_marginals}] We now
prove the statements concerning $Q_t$ and $Q_t^\eta$. We need
only to show that the measures, conditional on $\eta$, are equivalent since
the full measures are simply the integration of the conditioned
measures against the Wiener measure governing $\eta$. We will use
Girsanov's Theorem (cf.  \cite{b:Oksendal92,b:RevuzYor94} Ch 8, Thm
1.1) to compare \eqref{eq:toy_l} with the Ornstein-Uhlenbeck process
$\dt{z}=-\nu_2 z + \sigma_2 \dt{\xi}$.  Girsanov's Theorem states that
the two measures on path space are equivalent if a certain exponential
martingale, which gives the Radon-Nikodym derivative, is uniformly
integrable. This is guaranteed by Novikov's criterion (cf.
\cite{b:Oksendal92,b:RevuzYor94} Ch 8, Prop 1.15) which, translated
into our setting, becomes $\EE \exp\left( \frac12 \int_0^t
  \frac1{\sigma_2^2}|F_2(\ell,h)|^2 ds \right) < \infty$.  Since by
assumption
\begin{equation*}
\EE \exp\left( \frac12 \int_0^t
  \frac1{\sigma_2^2}\left|F_2\big(\ell(s),h(s)\big) \right|^2 ds
\right)\leq\exp\left( 
  \frac{K^2}{2\sigma_2^2}t \right) < \infty,  
\end{equation*}
we know that the measures induced on path space by $\ell_{[0,t]}$
conditioned on $\eta$ and $z_{[0,t]}$ are equivalent. This in turn
implies that the time $t$ marginals are equivalent. Since the law of
$z(t)$ for fixed $t$ is Gaussian and thus equivalent to Lebesgue
measure the proof is complete.  \epf

\bpf[Proof of the second part of Lemma \ref{l:toy_marginals}]  We now
prove the statements concerning $Q_{[0,\infty)}$ and $Q_{[0,\infty)}^\eta$. Again we
use Girsanov's Theorem and only consider the conditioned measures.
This time we compare the measures induced on $[0,t]$ by $\ell$
starting from the same $\ell_0$ with the same $\eta$ but different
$h$'s. In this case, Novikov's criterion becomes
\begin{multline*}
  \EE \exp\Bigl( \frac12 \int_0^t
  \frac1{\sigma_2^2}
        |F_2(\ell,\Phi_{0,s}^\eta(\ell_{[0,s]};h_0))-
  F_2(\ell,\Phi_{0,s}^\eta(\ell_{[0,s]};\tilde h_0))|^2 ds \Bigr)\\ 
        \leq \exp\Bigl( \frac12 \int_0^t \frac{L_2}{\sigma_2^2}
    |\Phi_{0,s}^\eta(\ell_{[0,s]};h_0) -
    \Phi_{0,s}^\eta(\ell_{[0,s]};\tilde h_0)|^2 ds \Bigr)\\
        \leq \exp\Bigl( \frac12 \int_0^t \frac{L_2}{\sigma_2^2} |h_0
    -\tilde h_0|^2 e^{-2(\nu_1-L_1)s} ds \Bigr) 
\\
        \leq \exp\Bigl(
    \frac{L_2}{4\sigma_2^2(\nu_1-L_1)} |h_0 -\tilde h_0|^2 \Bigr).
\end{multline*}
Since the bound is finite and uniformly bounded in $t$, we conclude that the
measures on path space are equivalent on the time interval
$[0,\infty)$.  \epf

\subsection{Basic Ergodicity}
\label{sec:generalErgodicity}
We now present some general theorems which we will use to prove the
ergodicity of the SNS equations and the toy model. Hopefully, the
assumptions will seem natural in light of the structure of the toy model.
\begin{assumption}\label{a:contraction} There exists a set
  $\mathcal{B}\subset C(0,\infty;\X)$, with $\PP\{
  \varphi_{[0,\infty)}^W u_0 \in \mathcal{B} \} =1$ for all $u_0 \in
  \X$ so that the following holds: 
  
  If $u(t,W)$ and $\tilde u(t,W)$ are solutions to
  \eqref{eq:abs_equation} and  $\Omega_0$ a subset of $\Omega \times
  \Omega$ so that 
  \begin{align*}
    (W,\tilde W) \in \Omega_0 \Longrightarrow 
    \begin{cases}
      u(\ccdot,W)  ,\tilde u(\ccdot,\tilde W) \in \mathcal{B}\\
       \Ph W(t)-\Ph W(0) = \Ph \tilde W(t)- \Ph
  \tilde W(0)  & \mbox{ for all } t\geq 0\\
  \Pl u(t,W) = \Pl \tilde
  u(t,\tilde W) & \mbox{ for all } t\geq 0
    \end{cases}
   \end{align*}
   then $\normX{\Ph u(t,W) -  \Ph \tilde u(t,\tilde W)}\rightarrow 0$
  as $t\rightarrow \infty$ for all  $(W,\tilde W) \in \Omega_0$.
\end{assumption}
In the toy model the set $\mathcal{B}$ was not needed; the conclusion
held for all paths. This is not true in more general settings; we need
to restrict ourselves to a set of ``nice'' paths. For the SNS
equation, $\mathcal{B}$ will be the set of paths which grow and
average in a typical fashion. Recall from \eqref{eq:defQ}, that
$Q_{t}^\eta(u_0,\ccdot)$ and $Q_{[t,\infty)}^\eta(u_0,\ccdot)$ are
respectively the measure induced on the ``low modes'' $\Xl$ at time
$t$ by $u(t)$ and on the path space $C([0,\infty),\Xl)$ by
$u_{[t,\infty)}$ if one conditions to use the noise realization $\eta$
and to start from the initial condition $u_0$ at time $t=0$.

\begin{assumption} \label{a:leaves} 
  For all $u_0=(\ell_0,h_0) \in\Xl\times\Xh$,
  $Q_t^\eta(\ell_0,h_0,\ccdot)$ is equivalent to Lebesgue measure for
  almost every $\eta$.  For all $u_0=(\ell_0,h_0)$ and $\tilde
  u_0=(\tilde \ell_0,\tilde h_0) \in\Xl\times\Xh$, the measure
  $Q_{[0,\infty)}^\eta(\ell_0,h_0,\ccdot)$ is equivalent to
  $Q_{[0,\infty)}^\eta(\ell_0,\tilde h_0,\ccdot)$ for almost every
  $\eta$.
\end{assumption}
As noted in the analysis of the toy problem, the equivalence of the
measure conditioned on $\eta$ implies the equivalence of the
unconditioned versions.

\begin{theorem}
\label{thm:ergodic}
If Assumptions \ref{a:contraction} and \ref{a:leaves} hold, then
\eqref{eq:abs_equation} has at most one $\X$-valued invariant
probability measure.
\end{theorem}

By an $\X$-valued probability measure $\mu$, we mean a measure such
that $\mu(\X)=1$.  Once this theorem is proven, we will have proven
the ergodicity of the toy problem from the previous section. In the
SNS setting, notice that Assumption \ref{a:leaves} is close to
Theorem \ref{thm:determiningModes}. Lemma \ref{l:abs_contraction}
makes the set $\mathcal{B}$ explicit.  We now state a number of lemma
which will be used to prove Theorem \ref{thm:ergodic}.

\begin{lemma} \label{l:Equiv}
  Assume Assumption \ref{a:leaves} holds.  For any pair of initial
  conditions $u_0=(\ell_0,h_0)$ and $\tilde u_0=(\tilde \ell_0,\tilde
  h_0) \in\Xl\times\Xh$ and any $t > 0$, the measure
  $Q_{[t,\infty)}(\ell_0,h_0,\ccdot)$ and $Q_{[t,\infty)}(\tilde
  \ell_0,\tilde h_0,\ccdot)$ are equivalent.  Similarly for almost
  every $\eta$, $Q_{[t,\infty)}^\eta(\ell_0,h_0,\ccdot)$ is equivalent
  to $Q_{[t,\infty)}^\eta(\tilde \ell_0,\tilde h_0,\ccdot)$.
\end{lemma}

Given any invariant measure $\mu$, we define two classes of associated
measures; one on the future trajectories and one on the past
trajectories.  Let $\mathcal{M}_-$ denote the natural measure on
$C((-\infty,0];\X)$ defined by cylinder sets of the type: for some
$t_0,t_1,\cdots t_n, t_0<t_1<t_2\cdots t_n\leq0$,
\begin{equation*}
  A=\left\{\big(\ell(s),h(s)\big)\in
    C\big((-\infty,0],\X\big),\big(\ell(t_i),h(t_i)\big) \in A_i,
    i=0,\cdots n\right\}
\end{equation*}
where the $A_i$'s are Borel sets of $\X$. The definition
\begin{equation*}
  \mathcal{M}_-(A)= \mu(A_0)\cdot \PP\{ u(t_i) \in A_i , i=1\dots n |  u(t_0) \in A_0\} .
\end{equation*}
characterizes the measure.  Similarly we define $\mathcal{M}_+$ on
$C([0,\infty);\X)$.  We also define $\mathcal{M}_+^\eta$ on
$C([0,\infty);\X)$ by pushing $\mu$ forward under dynamics conditioned
to use the noise realization $\eta$. We will define a measure
$\mathcal{M}_-^\eta$ at the end of the section and explore it
properties. Recalling that $\Pl$ was the projection onto $\Xl$; and
hence, $\mathcal{M}_+\Pl^{-1}$ is a measure on $C([0,\infty);\Xl)$.
Analogously for $t>s$, we define $\Pi_{[s,t)}$ as the projection onto
the space $C([s,t);\X)$. 
\begin{assumption} \label{a:equivFutures}
  Let $\mu_1$ and $\mu_2$ be any two invariant measures and let
  $\mathcal{M}_{+,1}^\eta$ and $\mathcal{M}_{+,2}^\eta$ be the
  measures induced on $C([0,\infty);\X)$ described above. Then for
  almost every $\eta$, $\mathcal{M}_{+,1}^\eta\Pl^{-1}$ is equivalent to
  $\mathcal{M}_{+,2}^\eta\Pl^{-1}$.
\end{assumption}

\begin{lemma}\label{l:stationaryPath}
  Assumption \ref{a:leaves} implies Assumption \ref{a:equivFutures}.
\end{lemma}
\bpf[Proof of Lemma \ref{l:stationaryPath}] Since $\mu$ is invariant,
for any $B \subset C([t,\infty);\Xl)$ and $t > 0$
\begin{align*}
\mathcal{M}_{+}^{\eta}\Pl^{-1}(
B) = \int_{\X} Q_{[t,\infty)}^{\theta_{-t}\eta}(\ell_0,h_0,B)\mu(d\ell_0\times dh_0)  
\end{align*}
the result follows from Lemma \ref{l:Equiv} and since $\theta_t$ is
ergodic; hence, mapping one set of full measure in $\Pl \Omega$ to another set of
full measure.\epf

Assumption \ref{a:equivFutures}, is weaker than Assumption
\ref{a:leaves}. As the next lemma shows, it is sufficient to prove
ergodicity. In some settings where solutions to the initial value
problem do not have nice moment properties it is more convenient to work
directly with stationary solutions. This type of analysis is presented
in \cite{b:BakhtinMattingly03Pre}. However, in systems like the SNS
equations such difficulties do not arise and, as we shall see,
Assumption \ref{a:leaves} holds.

In light of above lemma, the following result implies Theorem \ref{thm:ergodic}.
\begin{lemma}\label{l:ergodic}
  If Assumptions \ref{a:contraction} and \ref{a:equivFutures} hold, then
  \eqref{eq:abs_equation} has at most one $\X$-valued invariant
  probability measure.  
\end{lemma}

Note: If one was only interested in events which depended on the part
of the path in $\Xl$ then Assumption \ref{a:contraction} is not
needed.

\bpf[Proof of Lemma \ref{l:ergodic}] Since all invariant measures are
a linear combination of ergodic measures it is enough to show there is
a unique ergodic measure. Let $\mu_1$ and $\mu_2$ be two different
ergodic measures. Let $\mathcal{M}_{+,1}$, $\mathcal{M}_{+,1}^\eta$,
$\mathcal{M}_{+,2}$ and $\mathcal{M}_{+,2}^\eta$ be the associated
measures defined above.  Let $\phi:\Xl \times \Xh \rightarrow \RR$ be
a measurable test function bounded with $\sup\normX{\phi}\leq 1$ and
$\phi(\ell,\ccdot) \in \mbox{Lip}_1(\Xh)$ for all $\ell$. The norm
induced on measures by this class of test functions dominates the
Wasserstein ( or Kantorovich) distance for measures. Hence, this class
of test functions is rich enough so that if $\int \phi d\mu_1=\int
\phi d\mu_2$ for all such $\phi$ then $\mu_1=\mu_2$ (cf.
\cite{b:Dudley76} ).  Since $\mathcal{M}_{+,i}$ is invariant under the
flow induced on measures, the Birkoff ergodic theorem implies that
there exists sets $\bar A_i \subset C([0,\infty),\X)$ such that
$\mathcal{M}_{+,i}(\bar A_i)=1$ and for all $(\ell,h) \in \bar A_i$
\begin{align}
  \label{eq:conv1}
  \lim_{t\rightarrow\infty}\frac1t \int_0^{t} \phi( \ell(s),h(s) ) ds
  =\bar \phi_i= \int \phi(x,y)\mu_i(dx \times dy) .
\end{align}
Define $A_i=\bar A_i \cap \mathcal{B}$ where $\mathcal{B}$ is the set
from Assumption \ref{a:contraction}. Again remark that
$\mathcal{M}_{+,i}(A_i)=1$ since $\mathcal{B}$ has full measure. Since
$\mathcal{M}_{+,i}(\bar A_i)=\EE\mathcal{M}^\eta_{+,i}(\bar A_i) = 1$,
$\mathcal{M}_{+,i}^\eta(A_i)=1$ for almost every $\eta$.  Let
$A_i^\eta$ be a subset of $A_i$ of full
$\mathcal{M}_{+,i}^\eta$--measure so that the paths in $A_i^\eta$ are
solutions with a noise realization $W$ so $\Pl W=\eta$.  By Assumption
\ref{a:equivFutures}, $\mathcal{M}_{+,1}^\eta \Pl^{-1}$ is equivalent
to $\mathcal{M}_{+,2}^\eta \Pl^{-1}$. Hence, $\mathcal{M}_{+,2}^\eta
\Pl^{-1}(\Pl A_1^\eta)>0$ which implies that $\Pl A_1^\eta \cap \Pl A_2^\eta$ is not
empty since $\mathcal{M}_{+,2}^\eta\Pl^{-1}(\Pl A_2^\eta)=1$.  Hence, the
set $A=\{(\ell,h_1,h_2): (\ell,h_1) \in A_1^\eta \mbox{ and } (\ell,h_2)
\in A_2^\eta\}$ is not empty. Fixing some $(\ell,h_1,h_2) \in A$, from
\eqref{eq:conv1} and Assumption \ref{a:contraction} we have that for
any $\epsilon>0$ there exists a $T$ so that for all $t>T$
\begin{align*}
  \Bigl|\frac1t \int_0^{t} \phi( \ell(s),h_i(s) ) ds
  - \bar \phi_i\Bigr| < \frac\epsilon4
\end{align*}
and $\normX{h_1(t)-h_2(t)} <\frac\epsilon4$. This last inequality
holds because the hypotheses of Assumptions \ref{a:contraction} are satisfied.
Hence,
\begin{align*}
  |\bar \phi_1 - \bar \phi_2| &\leq \sum_i \Bigl| \frac{1}{t} \int_0^{t}
   \phi(\ell(s),h_i(s))ds - \bar \phi_i\Bigr| + \Bigl| \frac{1}{t} \int_0^{t}
   \phi(\ell(s),h_1(s)) -  \phi(\ell(s),h_2(s)) ds\Bigr|\\
 & \leq \frac\epsilon2 + \frac{2T}{t} + \frac\epsilon4 \frac{t-T}{t} <
   \epsilon \mbox{ for $t$ sufficiently large.}
\end{align*}
Since $\epsilon$ was arbitrary, the proof is complete.
\epf

\bpf[Proof of Lemma \ref{l:Equiv} ] 
Let $A \subset C(0,\infty;\X)$. We will show that 
\begin{gather*}
        Q_{[t,\infty)}^\eta(\ell_0,h_0,A)=0
                \text{ implies that }
                Q_{[t,\infty)}^\eta(\tilde\ell_0,\tilde h_0,A)=0.
\end{gather*}
 First notice that for $B \subset \Ph\X$ if
$H_t^\eta(\ell_0,h_0,x,B)=\PP\{ h(t) \in B | \ell(0)=\ell_0, h(0)=h_0,
\ell(t)=x, \mathcal{F}^\eta_{[0,t]}\}$ for $t>0$ then 
\begin{align}
  \label{eq:toy_equiv}
  Q_{[t,\infty)}^\eta(\ell_0,h_0,A) = \int \int
  Q_{t}^\eta(\ell_0,h_0,dx)H_t^\eta(\ell_0,h_0,x,dy)
  Q_{[0,\infty)}^{\theta_t\eta}(x,y,A) \ .
\end{align}
Hence for almost every $\eta$, $Q_{[t,\infty)}^\eta(\ell_0,h_0,A)=0$
implies that $Q_{[0,\infty)}^{\theta_t\eta}(x,y,A)=0$ for
$\mbox{Leb}(dx)\times H_t^\eta(\ell_0,h_0,x,dy)$ almost every $(x,y)$
because by assumption $Q_{t}^\eta(\ell_0,h_0,\ccdot )$ is equivalent
to Lebesgue measure. By the second part of Assumption \ref{a:leaves},
we know that $Q_{[0,\infty)}^{\tilde \eta}(x,y,\ccdot)$ is equivalent
to $Q_{[0,\infty)}^{\tilde \eta}(x,\tilde y,\ccdot)$ for all
$x,y,\tilde y$ and for $\PP$-almost every $\tilde \eta$. Hence,
$Q_{[0,\infty)}^{\theta_t \eta}(x,y,A)=0$ for $\mbox{Leb}(dx)\times
H_t^\eta(\tilde \ell_0,\tilde h_0,x,dy)$ almost every $(x,y)$ and
$\PP$-almost every $\eta$. (Here we have used that the shift is
ergodic with respect to $\PP$. So $\theta_t$ a maps set of full
measure to another set of full measure.)  And hence, by the
representation for $Q_{[t,\infty)}^\eta(\tilde \ell_0,\tilde h_0,A)$
analogous to \eqref{eq:toy_equiv}, we conclude
$Q_{[t,\infty)}^\eta(\tilde \ell_0,\tilde h_0,A)=0$. \epf

\subsection{One Force, One Solution: Statistical Equilibrium \\Measures
  and Trivial Random Attractors} In analogy to $\mathcal{M}^\eta_+$,
we define $\mathcal{M}^\eta_-$ on $C((-\infty,0];\X)$ as the limit as
$t\rightarrow - \infty$ of
$M_{t,0}^\eta=\EE\{\varphi_{[t,0]}^{(\xi,\eta)}\mu|\mathcal{F}_{(-\infty,0]}^\eta\}$.
As discussed in \cite{b:LeJan87,b:DarlingYves88,b:bax91}, the sequence
is a backwards martingale, hence the limit exists almost surely by the
martingale convergence theorem. By $\varphi_{[t,0]}^{(\xi,\eta)}u_0$
we mean the entire piece of trajectory on $[t,0]$.  Similarly, one can
define $\mathcal{M}^{(\xi,\eta)}_-$ on $C((-\infty,0];\X)$ by the
limit $t \rightarrow -\infty$ of $\varphi_{[t,0]}^{(\xi,\eta)}\mu$. This
is the called the equilibrium measure \cite{b:LeJan87} and $
\mathcal{M}^\eta_-= \EE\{\mathcal{M}^{(\xi,\eta)}_-|
\mathcal{F}_{(-\infty,0]}^\eta\}$.  In a similar manner, one can
define $\mathcal{M}^\eta$ on all of $C(-\infty,\infty; \Xl)$.

In the case of the SNS equations Flandolli, Craul, and Debussche
\cite{b:Fl94,b:CrauelDebusscheFlandoli97} proven the existence of a
compact random attractor $\mathcal{A}(\xi,\eta) \subset \L^2$ which
attracts all bounded subsets $B \subset \L^2$ in the sense
\begin{align*}
  \lim_{t\rightarrow -\infty}
  d(\varphi_{t,0}^{(\xi,\eta)}B,\mathcal{A}(\xi,\eta))=0
\end{align*}
where $d$ is the symmetric Hausdorff distance on sets.

If we define the action of the shift $\theta_t$ of $\eta$ as
$(\theta_t\eta)(s)=\eta(t+s)-\eta(t)$, the measure $\mathcal{M}^\eta$
is invariant under the skew flow on measures fibered over $\eta$.
Then we have, for example, $\EE\{ \varphi^{(\xi,\eta)}_{0,t}
\mathcal{M}^\eta | \mathcal{F}^\eta\} = \mathcal{M}^{\theta_t
  \eta}$. (Recall that $ \mathcal{F}^\eta$ was the sigma algebra
generated by $\eta$.)
Similarly, when a random attractor exists
$\varphi^{(\xi,\eta)}_{0,t}\mathcal{A}(\xi,\eta) =
\mathcal{A}(\theta_t\xi,\theta_t\eta)$.

We can consider the equation \eqref{eq:abs_h} in isolation over a
probability space $\Ph\Omega \times C(-\infty,\infty;\Xl)$ with the
measure $\PP_\eta(d\eta) \times \mathcal{M}^\eta\Pi_\ell^{-1}(d\ell)$.
In other words, we have elevated the part of the phase space
$C(-\infty,\infty;\Xl)$ to part of the base probability space. On this
space the $h(t)$ dynamics has the same property as the whole SNS
equation under the extremely contractive assumption( $N_*^2 \geq
\mathcal{C}\frac{\mathcal{E}_0}{\nu^3}$). In particular, an analogous
theorem to Theorem \ref{thm:largeNu} holds: there is a unique solution
$h_*(t;\eta,\ell)$ which attracts all other solutions.  In these
coordinates, the random attractor for the equation \eqref{eq:abs_h} is
the single solution $h_*(t;\eta,\ell)$. Therefore, the invariant
measure $\mu$ from above projected onto $\Xh$ disintegrates into a
delta measure concentrated at $\delta_{h_*(0;\eta,\ell)}$ against the
measure $\PP_\eta(d\eta) \times
\mathcal{M}^\eta\Pi_\ell^{-1}(d\ell)$. That is to say, if $F:\Xl
\times \Xh \rightarrow \RR$ then
\begin{align*}
  \int F(\ell,h)\mu(d\ell \times dh) &= \int F(\ell,h_*(0;\eta,\ell))\big[
\PP_\eta(d\eta) \times\mathcal{M}^\eta\Pi_\ell^{-1}(d\ell)\big]\\ 
&= \int F(\ell,h)\delta_{h_*(0;\eta,\ell)}(dh)\big[ \PP_\eta(d\eta)
\times\mathcal{M}^\eta\Pi_\ell^{-1}(d\ell)\big] .
\end{align*}
This is the analog for the partially dissipative system of the ``one
force, one solution'' (i.e. trivial random attractor) discussed in
\cite{b:Mattingly98,b:EKhaninMazelSinai00,
  b:Schmalfuss97,b:LeJan87,b:Mattingly02c,b:MasmoudiYoung02,b:EVandenEijnden00a}
or exemplified by Theorem \ref{thm:largeNu}. A similar statement holds
for the toy problem and all of the systems satisfying the assumptions
in section \ref{sec:generalErgodicity}.

\section{Contractive Nature of the SNS Dynamics}
\label{sec:contractive}

The proof of ergodicity of the SNS under the assumption that only the
``determining modes'' are forced will parallel the proof of the toy
model. Our first step is to establish Assumption \ref{a:contraction}
in the context of the SNS. We do this by proving a quantitative
version of Theorem \ref{thm:determiningModes} which was given earlier.

To see what is involved we consider two solutions to equation
\eqref{eq:highMode}, $h_1$ and $h_2$, driven by a common low mode
process $\ell$ and noise $\eta$. That is for $t>s$,
$h_i(t)=\Phi_{s,t}^\eta(\ell_{[s,t)};h_i(s))$. Denoting
$u_i=\ell+h_i$, we have
\begin{align*}
  \pdt{\rho(t)}=&-\Lambda^2 \rho(t) + \Ph B(u_1,\rho) + \Ph
  B(\rho,u_2) \intertext{which using standard estimates on the
    nonlinearity (cf.  \cite{b:CoFo88}) produces, for some
    $\mathcal{C}>0$, } \dt{\LL{\rho}^2}\leq& -2\nu\LL{\Lambda \rho}^2
  + (2\mathcal{C})^\frac12\LL{\Lambda \rho}\LL{\Lambda u_2}\LL{\rho}
  \leq -\nu\LL{\Lambda
    \rho}^2 + \frac{\mathcal{C}}{\nu} \LL{\Lambda u_2}^2\LL{\rho}^2\\
  \leq& -\left(\nu N_*^2 - \frac{\mathcal{C}}{\nu} \LL{\Lambda
      u_2}^2\right)\LL{\rho}^2.
\end{align*}
The above estimate then
gives for $t_0< t$
\begin{align}\label{eq:SNSContractive}
\LL{\rho(t)}^2 \leq &\LL{\rho(t_0)}^2 \exp\left( -\nu N_*^2(t-t_0) +
  \frac{\mathcal{C}}{\nu}\int_{t_0}^t \LL{\Lambda u_2(s)}^2ds
  \right)
\end{align}
We now see the new difficulty which the Stochastic Navier Stokes
equations present over the toy model. The contraction rate
depends on the time average of the enstropy $\LL{\Lambda u_2(s)}^2$ of
one of the solutions. However, after we develop some estimates
controlling this quantity the proof will proceed using standard ideas
of localization from stochastic analysis.

\section{The Energy and Enstrophy}
\label{sec:local}

The toy model is an extremely uniform setting. The added difficulty
in the SNS relative to the toy model, is the lack of uniformity.
However, the standard idea of localization from stochastic analysis
allows us to overcome this hurdle. As we saw in the last section, the
growth of the energy $\LL{u}^2$ and the time average of the enstrophy
$\LL{\Lambda u}^2$ seem to be of importance in controlling the
uniformity of the contraction. This will be come clearer after the
next two sections. We begin with some estimates on the energy and enstrophy.
\begin{lemma}
  \label{l:uniformEnergyEnstrophy}
  $\EE \LL{u(t)}^{2} \leq e^{-2 \nu(t-t_0)} \EE \LL{u(t_0)}^{2} +
  \frac{\mathcal{E}_0}{2\nu}\left(1- e^{-2 \nu ( t-t_0)}\right)$ and
  for any $p \geq  1$, $\EE \LL{u(t)}^{2p} \leq e^{-2 \nu (t-t_0)} \EE
  \LL{u(t_0)}^{2p} + C_0 \int_{t_0}^t e^{-2 \nu (t-s)} \EE
  \LL{u(s)}^{2(p-1)} ds$.
\end{lemma}
This implies that if one has a solution $u(t,W)$ defined for $t \in
(-\infty,\infty)$ such that $e^{-2 \nu |t|}\EE
\LL{u(t,W)}^2\rightarrow 0$ as $t \rightarrow -\infty$ then in fact
$\EE \LL{u(t,W)}^2$ is uniformly bounded in time. Using similar
reasoning, one can show the following result.
\begin{lemma}
  \label{l:invMeasureFiniteEnergy}Assume that $\mu$ is an invariant
  measure such that there exist a $U \subset \L^2$ with $\mu(U)=1$.  For
  any such measure stationary measure all energy moments are finite. In
  fact for any $p\geq 1$ there exist constant $C_p< \infty$ such that
  $\int_{\Ltwo} \LL{u}^{2p}d\mu(u) < C_p$ for all invariant measures
  $\mu$. In particular, $C_1=\frac{\mathcal{E}_0}{2\nu}$. Furthermore,
  $\int_{\L^2} \LL{\Lambda u}^2 d\mu(u)= \frac{\mathcal{E}_0}{2\nu}$
  assuming only that $\mathcal{E}_0=\sum |\sigma_k|^2 < \infty$. If
  $\mathcal{E}_1=\sum|k|^2 |\sigma_k|^2 < \infty$ is finite then the
  analogous statements hold for $ \LL{\Lambda u}^{2p}$ replacing
  $\LL{u}^{2p}$. In particular, $ \int_{\L^2} \LL{\Lambda^2 u}^2
  d\mu(u)= \frac{\mathcal{E}_1}{2\nu}$.
\end{lemma}
Since one can construct a stationary solution from any invariant
measure and vice versa (see section \ref{sec:generalErgodicity} ),
this conclusion applies equally to any stationary solution.  The
proofs of Lemma \ref{l:uniformEnergyEnstrophy} and
\ref{l:invMeasureFiniteEnergy} can be found in the appendix of
\cite{b:EMattinglySinai00}. Related statements can be found in Chapter
3 and 4 of \cite{b:Mattingly98b}, section 2 of \cite{b:Mattingly98},
or the appendix of \cite{b:Mattingly02a}. The moment estimates are
just the stochastic analogs of deterministic estimates. Similar
estimates from slightly different points of view can be found in
\cite{b:VishikFursikov88,b:MikuleviciusRozovskii02,b:BricmontKupiainenLefevere00}.
If one assumes, $\mathcal{E}_1=\sum|\sigma_k|^2|k|^2 < \infty$ then
completely analogous statements can be made about the enstrophy.

It is critical to our analysis to understand the typical size of the
fluctuations of the enstrophy about its mean of
$\frac{\mathcal{E}_0}{2\nu}$. Applying It\^o's formula to the energy,
one obtains
\begin{align*}
  d \LL{u(s)}^2 = - 2\nu\LL{\Lambda u(s)}^2 dt + \mathcal{E}_0 dt +
 2 \ip{u}{dW(t)}_{\L^2}
\end{align*}
If one writes the last term as
\begin{align*}
  \left(\sum_k |\sigma_ku_k|^2\right)^\frac12
  \left[\frac{\ip{u}{dW(t)}_{\L^2}}{ \left(\sum_k |\sigma_ku_k|^2\right)^\frac12}\right],
\end{align*}
then the term in the square brackets is distributed as a one
dimensional Brownian motion adapted to the filtration generated by the
$W$ increments. This motivates our definition of a Lyapunov function
in the abstract setting \eqref{eq:abs_equation}, which is contained in
the next section.

\section{Growth and Fluctuations in A General Setting}
\label{sec:fluctuationsGeneral}
In this section, we put an abstract framework on the ideas of
the previous section. In the section \ref{sec:highViscosity}, we return to the concrete
setting of the SNS.
\begin{assumption}\label{a:Lyop}
There exists a function $V:\X \to [0, \infty)$, with $c_0 V(x)^{p_0}
   \geq \normX{x}^2$ for positive $p_0,c_0$, so that
   for a solution $u$ of equation \eqref{eq:abs_equation} $V$
    satisfies the It\^o equation
    \begin{equation*}
      dV(u(t)) = g(u(t))dt + f(u(t)) d\widetilde B(t).
    \end{equation*}
    Here $\widetilde B$ is a standard one dimensional Wiener process
    adapted to the flow generated by $(dW)$.  $g:\X\to\RR\cup\infty$ is a
    function satisfying
    \begin{equation*}
      g(u)<C_1 - U(u) \mbox{ where }   U(u) \geq C_2 V(u)
    \end{equation*}
    for some constants $C_1,C_2 >0$ and $U:\X \to [0,\infty)\cup\infty
    $.  Though $U$ and $g$ might be infinite on $\X$, we assume if
    $u(t)$ is a solution to \eqref{eq:abs_equation} on $[0,t]$, then
    $\int_0^tU(u(s)) ds < \infty$ almost surely and the above
    inequalities holds whenever $U(u) < \infty$.  And $f:\X \to
    \RR\cup\infty$ is a function satisfying\label{er:GV1}
    \begin{equation*}
       C_3|f(u)|^2 \leq U(u)
    \end{equation*}
    for some $C_3>0$.
     \end{assumption}
 From the calculation at the end of the last section, in the SNS
 setting we should take $V(u)=\LL{u}^2$, $U(u)=2\nu\LL{\Lambda u}^2$,
 $C_1=\mathcal{E}_0$, $C_2=2\nu$, and $C_3=\frac{\nu}{2\sigma_*^2}$.
 (We could have also used $V(u)=\LL{\Lambda u}^2$,
 $U(u)=2\nu\LL{\Lambda^2 u}^2$ if $\mathcal{E}_1=\sum |\sigma_k|^2 |k|^2 < \infty$.)\label{er:E2}

\begin{lemma} \label{l:growthForward} For any $\epsilon \in (0,1)$ and
  $K>
\frac{1}{\epsilon C_3}$ one has 
  \begin{align*}
    \PP\left\{ \sup_{t>0}  \frac{ V(u(t)) + (1-\epsilon)\int_0^t
    U(u(s)) ds - V(u_0) - C_1 t}{1+\log(1+t)} > K
    \right\} \leq \exp\Bigl(-2 C_3 \epsilon K\Bigr)  
  \end{align*}
\end{lemma}
\bpf
Let $M(t)$ denote the martingale $\int_0^t f(u(s)) d\widetilde
W(s)$. Its quadratic variation $[M,M](t)$ is $\int_0^t f(u(s))^2
ds$. Since 
$C_3 [M,M](t) \leq \int_0^t U( u(s)) ds$, by It\^o's formula we have
\begin{align*}
  V(u(t)) + (1-\epsilon)\int_0^t U(u(s))ds - V(u_0) - C_1 t \leq M(t)-
    \epsilon C_3  [M,M](t)
\end{align*}
The exponential martingale estimate implies that 
\begin{equation}\label{eq:expMartEstimate}
  \PP\Big\{ \sup_{t \leq T}  M(t)-
    \epsilon C_3 [M,M](t) > a \Big\} \leq e^{- 2 \epsilon
    C_3 a } .
\end{equation}
Setting $a=K[1+\log(T)]$ one sees that probability of the event in the
statement of the lemma is bounded from above by
\begin{gather*}
\sum_{n=1}^{\infty}  \exp\bigl(-2 C_3 \epsilon
  K[1+\log(n)]\bigr) \leq \int_1^\infty \exp\bigl(-2 \epsilon C_3
  K-2\log(x)\bigr)dx = \exp\bigl(-2 C_3 \epsilon K\bigr)  .
\end{gather*}

More details can be found in the proofs of the following related
results: Theorem 4--6 in \cite{b:BakhtinMattingly03Pre}, Lemma A.5
and Lemma B.3 of \cite{b:EMattinglySinai00}, or for the use of the
exponential martingale Lemma A.2 \cite{b:Mattingly02}.
\epf

Using similar reasoning one can prove (see for example
\cite{b:BricmontKupiainenLefevere00,b:Mattingly02a}):
\begin{lemma}
  There exist positive constants $\gamma$ and $K$ so that for all invariant
  measures $\mu$ with $V(u) < \infty$ $\mu$-almost surely, $\int
  \exp(\gamma V(u)) d\mu(u) \leq K < \infty$ and for every initial
  condition $u_0$ and $t\geq 0$, $\EE \exp(\gamma V(u(t))) \leq K \exp(\gamma V(u_0))$
\end{lemma}
We now give estimates backward in time for stationary solutions.
\begin{lemma} \label{l:growthBackward}
  Let $u(t,W)$ be a stationary solution to \eqref{eq:abs_equation}
  with $V(u(t)) < \infty$ almost surely. There exists a $K_0$ and a
  $\gamma >0$ so that for $K > K_0$ 
  \begin{align*}
    \PP\Bigg\{ \sup_{t\in \RR} \frac{V(u(t)) +
      (1-\epsilon)\left|\disp\int_0^t U(u(t)) ds\right| -
      C_1|t|}{1+\log(1+|t|)} > K \Bigg\} \leq \exp(-\gamma K)
  \end{align*}
\end{lemma}
\bpf
The proof is essentially the same as that of Lemma \eqref{l:growthForward}. We
write
\begin{multline*}
  \sup_{t \in [T-1,T]} V(u(t)) + (1-\epsilon)\left|\disp\int_0^t
      U(u(t)) ds\right| - C_1|t| \leq V(u(T-1)) +\\
    \sup_{t \in [T-1,T]}\left[V(u(t)) 
      + (1-\epsilon)\left|\int_0^t U(u(s))ds\right| -
      C_1|t| -V(u(T-1))\right]
  \end{multline*}
By the previous two lemmas both of these terms have exponential
moments uniform in $T$. Using the same reasoning as in the end of the
proof of Lemma \ref{l:growthForward} completes to proof.\epf

In light of Lemmas \ref{l:growthForward} and \ref{l:growthBackward} we
define the following sets of ``nice'' trajectories which average well
and grow in a typical fashion. Fixing some fixed $\epsilon_* \in
(0,1)$, which will be set differently in different contexts, we define
\begin{multline}\label{eq:A_nB_n}
  A_n
=\Bigl\{  u\in C(-\infty,0; \X): \sup_{t \leq 0}
  \frac{V(u(t)) + (1-\epsilon_*)\bigl|
\int_0^t
      U( u(s)) ds\bigr| - C_1 t}{1+\log(t+1)}
  \leq n \Bigr\},\\
\shoveleft{
  B_n
=\Bigl\{ u\in C(0,\infty; \X): 
}\\
\sup_{t \geq 0}
  \frac{V(u(t)) +(1-\epsilon_*)\bigl|
\int_0^t U(u(s))ds\bigr| -
    V(u_0) - C_1t}{1+\log(1+t)} \leq n \Bigr\}.
\end{multline}
The previous lemmas imply that, with probability one, any stationary
solution is contained in $\cup  A_n$ and the solution to
any initial value problem is contained in $\cup  B_n$. 

From these lemmas it is clear that
\begin{equation*}
  \lim_{t \rightarrow
    \infty}\frac1t\int_0^t U(u(s))ds \leq C_1\qquad\mbox{ and }\qquad\lim_{t \rightarrow
    -\infty}\frac1{|t|}\int^0_t U(u(s))ds \leq C_1   
\end{equation*}
almost surely. This is the result analogous to the final conclusion
of Lemma \ref{l:invMeasureFiniteEnergy}. But notice that Lemma
\ref{l:growthForward} and  \ref{l:growthBackward}
also give information about the size of the fluctuations.

\section{Strongly Contractive Case SNS: Proof of Theorem \ref{thm:largeNu}}
\label{sec:highViscosity}

In the next section, we will take a more abstract point of view on the
contractive nature of the SNS equation and other SPDE. However,
first for illustrative reasons, we continue with the explicit
calculations began in section \ref{sec:contractive} and use them to
prove Theorem  \ref{thm:largeNu}.

\bpf[Proof of Theorem \ref{thm:largeNu}]
We begin by proving uniqueness. Let us assume that there are two
solutions $u(t,W)$ and $u_*(t,W)$ defined for all $t \in
(-\infty,\infty)$. Let $\rho(t,W)=u(t,W)-u_*(t,W)$. Both are governed
by \eqref{eq:NS}, so subtracting produces
\begin{align*}
  \pdt{\rho(t)}=&-\Lambda^2 \rho(t) + B(u,\rho) + B(\rho,u_*) .
  \intertext{Using the same estimates on the nonlinearity
    as in \eqref{eq:SNSContractive}, we obtain} 
\dt{\LL{\rho}^2}\leq&
  -\left(\nu N_*^2 - \frac{\mathcal{C}}{\nu} \LL{\Lambda
      u_*}^2\right)\LL{\rho}^2.
\end{align*}
Notice that this is the same estimate as was obtained in
\eqref{eq:SNSContractive}.  Since $u_*$ is a solution, we know that it
is contained in some $A_n$, where $A_n$ was the set of ``nicely''
fluctuating and growing paths defined in the last section. Recall that
for the SNS: $V(u)=\LL{u}^2$, $U(u)=2\nu\LL{\Lambda u}^2$,
$C_1=\mathcal{E}_0$. Hence, $u \in A_n$ implies that for $t\leq 0$
\begin{align*}
    \LL{u(t)}^2 + (1-\epsilon_*)2\nu\left|\disp\int_0^t
      \LL{\Lambda u(s)}^2 ds\right| - \mathcal{E}_0| t| < n[1+\log(|t|+1)]
\end{align*}
where $\epsilon_*\in(1,0)$ is a free parameter which we will set momentarily.
Continuing the estimation of $\rho$, using this bound, produces for $t_0 <0$
\begin{align*}
  \LL{\rho(0)}^2 \leq &\LL{\rho(t_0)}^2 \exp\left( -\nu N_*^2|t_0| +
    \frac{\mathcal{C}}{\nu}\int_{t_0}^0 \LL{\Lambda u_2(s)}^2ds
  \right)\\
  \leq &\LL{\rho(t_0)}^2 \exp\left( -\Big[\nu N_*^2-
    \mathcal{C}\frac{\mathcal{E}_0}{\nu^2(1-\epsilon_*)}\Big]|t_0| +
    \frac{n[1+\log(|t_0|+1)]}{2\nu(1-\epsilon_*) } \right) .
\end{align*}
Picking $\epsilon_*$ so $\nu N_*^2-
\mathcal{C}\frac{\mathcal{E}_0}{\nu^2(1-\epsilon_*)}=\frac12\gamma_*$
where $\gamma_*=\nu N_*^2- \mathcal{C}\frac{\mathcal{E}_0}{\nu^2}$ and
using the assumption that the solutions are in $A_n$ to control
$\LL{\rho(t_0)}^2$ yields the estimate
\begin{align*}
  \LL{\rho(0)}^2 \leq 4[\mathcal{E}_0+ n(1+\log(1+|t_0|)] \exp\left(
    -\frac12\gamma_*|t_0| +
    \frac{n[1+\log(|t_0|+1)]}{2\nu(1-\epsilon_*) } \right).
\end{align*}
Taking $t_0 \rightarrow -\infty$, proves uniqueness. A similar
estimate shows that the solution to any initial value problem
converges exponentially forward in time to $u_0$.  The existence can
be deduced from the existence of a stationary measure; however, it is
instructive to construct it directly, which we now do.

Let $u_n(t)$ be the solution starting from initial value zero at time
$-n$. From Lemma \ref{l:growthForward}, we know that $\theta_{-n}u_n
\in B_{k_n}$ for some $k_n$. ($\theta_{-n}$ just shifts the path on
$[-n,\infty)$ to a path on $[0,\infty)$.) In addition, we know that
$\PP\{ k_n > n^\frac18 \} \leq \exp( - \frac{\nu
  \epsilon_*}{\sigma_*^2} n^\frac18)$. Hence, by the Borel-Cantelli
Lemma, there exists an $n_*$ so that $k_n \leq n^\frac18$ for all
$n> n_*$. Let $n>m> M > n_*$, then for $M$ sufficiently
large we have
\begin{multline*}
  \sup_{s\in [-1,0]}\LL{u_n(s) - u_m(s)} 
\leq \sum_{j=m+1}^n \sup_{s
    \in [-1,0]} \LL{u_{j+1}(s) - u_j(s)} 
\\\leq \sum_{j=M}^\infty
  \sup_{s \in
    [-1,0]} \LL{u_{j+1}(s) - u_j(s)}\\
 \leq \sum_{j=M}^\infty \LL{u_{j+1}(-j)} \exp(-\frac12
  \gamma_*|j-1| + |j|^\frac18[1 +
  \log(1+|j|)]) \\
 \leq \sum_{j=M}^\infty [\mathcal{E}_0 + 2|j|^\frac18]
  \exp(-\frac12 \gamma_*|j-1| + |j|^\frac18[1 + \log(1+|j|)])
\end{multline*}
Since the last sum is less than $C\exp(-\gamma M)$ for some positive
constants $C$ and $\gamma$, the sequence is Cauchy and the
proof is complete.
\epf

\section{Contractive Nature in the General Setting}\label{sec;genContraction}
We now extract the essential assumptions of the previous section and
present them in an abstract form. The choice of assumptions follows
\cite{b:BakhtinMattingly03Pre} which uses ever so slightly different
assumptions, but proves more detailed estimates.  In particular,
statements about the continuity of the map $\Phi$ are made. (See
Theorem 9 of \cite{b:BakhtinMattingly03Pre}.)  The treatment is also
informed and influenced by \cite{b:Hairer02,b:Mattingly02,b:ELui02}.
 
\begin{assumption} \label{a:contractive}Consider the $G:\X \to \X$ from
  \eqref{eq:abs_equation}. Assume that Assumption \ref{a:Lyop} holds
  and that for all $\ell \in \Xl$ and $h,\tilde h \in \Xh$, with
  $\ell+h,\ell+\tilde h \in \mathcal{D}(G)$,\label{er:G2} and some $c_i
  \geq 0$, and $p_i \geq 0$, with $c_1 > C_1(c_2+c_3)$,
  \begin{align*}
    \ipX{G(\ell + h)-G(\ell + \tilde h)}{h-\tilde h}&\leq
    \big[-c_1+c_2U(\ell+h)+c_3U(\ell+\tilde h)\big]\normX{h-\tilde h}^2\\
    \normX{\Pl G(\ell + h)-\Pl G(\ell + \tilde h)}^2&\leq c_4\left[
      1 + V(\ell + h)^{p_1}+ V(\ell + \tilde
      h)^{p_1}\right]\normX{h-\tilde h}^{p_2} \ .
  \end{align*}
\end{assumption}
  
We give the analog of Theorem \ref{thm:determiningModes} and Lemma
\ref{l:toy_contraction} in the general setting of equation
\eqref{eq:abs_equation}. This is a quantitative version of the
determining mode result given in Theorem \ref{thm:determiningModes}
and will be used to verify Assumption \ref{a:contraction}.
\begin{lemma}\label{l:abs_contraction}
  Let Assumption \ref{a:contractive} hold. In particular,
  $\gamma_*\eqdef c_1-C_1(c_2+c_3) >0$. Set the $\epsilon_*$ from the
  definition of $A_n$ and $B_n$ in \eqref{eq:A_nB_n}, so that
  $c_1-\frac{C_1(c_2+c_3)}{1-\epsilon_*} = \frac12 \gamma_*$ and
  define $c_*=\frac{c_2+c_3}{1-\epsilon_*}$.
  \begin{enumerate}
  \item Fixing a $T \in (0,\infty]$, let $\ell \in C(0,T;\Xl)$ and
    $h_0$, $\tilde h_0 \in \Xh$ and define $u(t)=\ell(t)+
    \Phi_{0,t}^\eta(\ell_{[0,t]};h_0)$ and $\tilde u(t)=\ell(t)+
    \Phi_{0,t}^\eta(\ell_{[0,t]};\tilde h_0).$

    Assume that  $u \in \Pi_{[0,T)} B_n$
    if $c_2>0$ and $\tilde u \in \Pi_{[0,T)} B_n$
    if $c_3>0$. Then for all $t \in [0,T)$
  \begin{align*}
    \normX{ \Phi_{0,t}^\eta(\ell_{[0,t]};h_0) -
   \Phi_{0,t}^\eta(\ell_{[0,t]};\tilde h_0)}^2 \leq\normX{h_0 -
   \tilde h_0}^2 e^{nc_*} (1+t)^{nc_*} e^{-\frac12\gamma_* t} \ .
  \end{align*}
\item 
   Fixing a $T \in (-\infty,0)$, let $\ell \in C(T,0;\Xl)$ and
    $h_0$, $\tilde h_0 \in \Xh$ and define $u(t)=\ell(t)+
    \Phi_{T,t}^\eta(\ell_{[T,t]};h_0)$ and $\tilde u(t)=\ell(t)+
    \Phi_{T,t}^\eta(\ell_{[T,t]};\tilde h_0).$
    
    Assume that $u \in \Pi_{[T,0]} A_n$ if $c_2>0$ and $\tilde u \in
    \Pi_{[T,0]} A_n$ if $c_3>0$. Then for all $t \in [T,0]$ 
      \begin{align*}
        \normX{ \Phi_{T,t}^\eta(\ell_{[T,t]};h_0) -
          \Phi_{T,t}^\eta(\ell_{[T,t]};\tilde h_0)}^2 \leq\normX{h_0 -
          \tilde h_0}^2 e^{nc_*} (1+|T|)^{nc_*} e^{-\frac12\gamma_*(|T|-|t|)} \ .
      \end{align*}
    \end{enumerate}
  \end{lemma}
This theorem can be restated in terms of solutions to \eqref{eq:abs_equation}.
\begin{corolary} \label{c:absContraction}Let $\Omega_0 \subset \Omega
  \times \Omega$ and $c_*$ and $\gamma_*$ are as defined in Lemma \ref{l:abs_contraction}.
  \begin{enumerate}
  \item Let $u(t,W)$ and $\tilde u(t,W)$ be solutions to
    \eqref{eq:abs_equation} on $C(0,T;\X)$ with $T> 0$.  If
    for $(W,\tilde W) \in \Omega_0$, one has $u(W), \tilde u(\tilde W)
    \in \Pi_{[0,T]}B_n$ and $\Ph W(s)-\Ph W(T)= \Ph \tilde W(s)-\Ph W(T)$ ,
    $\Pl u(s,W)=\Pl \tilde u(s,\tilde W)$ for $s \in [0,T]$ then
   \begin{equation*}
     \normX{ u(t,W)-\tilde u(t,\tilde W)}^2 \leq \normX{ u(0,W)-\tilde
     u(0,\tilde W)}^2  e^{nc_*}(1+t)^{nc_*} e^{-\frac12\gamma_* t} \ .
      \end{equation*}

  \item Let $u(t,W)$ and $\tilde u(t, W)$ be solutions to
    \eqref{eq:abs_equation} on
    $C(T,0;\X)$ with $T<0$.  If for $(W,\tilde W) \in
    \Omega_0$, one has $u(W), \tilde u(\tilde W) \in \Pi_{[T,0]}A_n$
    and $\Ph W(s)-W(T)= \Ph \tilde W(s)-W(T)$ , $\Pl u(s,W)=\Pl \tilde
    u(s,\tilde W)$ for $s \in [T,0]$ then
   \begin{multline*}
        \normX{ u(t,W)-\tilde u(t,\tilde W)}^2 \leq 2[C_1|T| +n+n\log(1+|T|)
          ]e^{nc_*} \times\\
(1+|T|)^{nc_*} e^{-\frac12\gamma_*(|T|-|t|)} \ .
      \end{multline*}

  \end{enumerate}
\end{corolary}

\bpf[Proof of Lemma \ref{l:abs_contraction} and Corollary
\ref{c:absContraction}] The proof of the two statements is almost
identical and is simply an abstraction of the ideas in the proof of
Theorem \ref{thm:largeNu} given in the last section. We give the
details of the first statement.

Let $\rho(s)= \Phi_{0,s}^\eta(\ell_{[0,s]};h_0) -
\Phi_{0,s}^\eta(\ell_{[0,s]};\tilde h_0)$, then equation
\eqref{eq:abs_h} and the assumption in the lemma and Assumption
\ref{a:contractive} imply that
\begin{align*}
  \frac12 \dt{\normX{\rho(s)}^2} & = \ipX{G(u)-G(\tilde
    u)}{\rho(s)}\\
  & \leq \big[-c_1 + c_2U(u(s)) + c_3U(\tilde u(s))
    \big]\normX{\rho(s)}^2  \ .
  \end{align*}
Since $u, \tilde u \in \Pi_{[0,T]}B_n$, we have
\begin{align*}
  \normX{\rho(t)}^2 &\leq
  \normX{\rho(0)}^2\exp\Bigl(-[c_1-\frac{C_1(c_2+c_3)}{1-\epsilon_*}]t
  +n\frac{c_2+c_3}{1-\epsilon_*}(1+\log(1+t))\Bigr)
\end{align*}
which proves the first result. The second result is just the same
except that the estimates from $A_n$ are used. See the proof of
Theorem \ref{thm:largeNu}. Corollary \ref{c:absContraction} is just a
restatement of the theorem with the added observation that backwards in
time the initial conditions $h_0=\Ph u(T,W)$ and $\tilde h_0=\Ph
\tilde u(T,\tilde W)$ can not grow too fast since the solutions are in
$A_n$.  \epf

Using the contractive properties backward in time one can define the
limit 
\begin{equation*}
  \lim_{t \rightarrow -\infty} \Phi_{t,0}^\eta(\ell;h_0)
\end{equation*}
for any $\ell$ which is a projection of a solution $u(t,W)$ on
$(-\infty,0]$. This limiting function, denoted $\Phi(\ell)$, is
independent of $h_0$ and can be used to reduce the dynamics to one on
$\Xl$ with memory (i.e. Gibbsian dynamics).  See
\cite{b:BakhtinMattingly03Pre} discussion of this in a general setting
and \cite{b:EMattinglySinai00,b:ELui02} for specific examples.  If one
endows $C(-\infty,0;\X)$ with the metric
$|u|_r=\sup_{t<0}\frac{\normX{u(t)}}{1+|t|^r}$, then in many settings
$\Phi(\ell)$ is continuous on the set of solutions $u$ with $\ell=\Pl
u$. In fact under some simple assumptions, it is globally Lipschitz on
each $B_n$ defined in section \ref{sec:fluctuationsGeneral}. In
particular, both of these facts hold for the SNS equation. See
\cite{b:BakhtinMattingly03Pre} for more discussion of this.

\section{Ergodicity: the SNS and the General Setting}
\label{sec:ergodicSNS}

We now turn to completing the proof of Theorem \ref{thm:mainErgodic}.
All that remains to prove is the last part of Theorem
\ref{thm:mainErgodic} about the ``essentially elliptic'' dynamics (the
case when $N_*^2 > \mathcal{C} \frac{\mathcal{E}_0}{\nu^3}$ ).  We
will do so by proving an ergodic theorem in the general setting of
\eqref{eq:abs_equation} and using the assumptions already introduced.

To prove basic ergodicity in this case, we will use Theorem
\ref{thm:ergodic}, which along with Corollary \ref{c:absContraction},
contains the essential ideas from \cite{b:EMattinglySinai00}. Lemma
\ref{l:growthForward} implies that almost every solution
$\varphi_{[0,\infty)}^W u_0$ is contained in a $B_n$ for some $n$.
This, coupled with Lemma \ref{c:absContraction}, is more than enough
to imply Assumption \ref{a:contraction} of Section
\ref{sec:generalErgodicity} with $\mathcal{B}=\cup B_n$.  We need only
verify Assumption \ref{a:leaves}, to prove complete the proof. Since
the author feels that techniques often used to verify the first part of
Assumption \ref{a:leaves} are suboptimal, we leave it as an assumption for the
moment. We will revisit the question at the end of this
section. Hence, we introduce the following assumption.
\begin{assumption}\label{a:badAssumptions}
  For all $t>0$ and $(\ell_0,h_0)\in \X$ and almost every $\eta$,
  $Q_t^\eta(\ell_0,h_0,\ccdot)$ is equivalent to Lebesgue measure.
\end{assumption}
The idea to prove the second part of Assumption \ref{a:leaves} is
again the idea of localization. By restricting ourselves to well behaved paths, we will be
able to of obtain the needed result for a subset of the probability
space. By relaxing the restriction, we can include arbitrarily large
subsets of the probability space, implying that the conclusion holds
with probability one. We prove the following result.
\begin{theorem} \label{thm:abs_ergodic}
  Consider equation \eqref{eq:abs_equation}. Let Assumption
  \ref{a:Lyop}, \ref{a:contractive} and \ref{a:badAssumptions} hold.
  In addition if $\sigma_k >0$ for all $k$ with $|k| \in (0,N_*)$,
  where $N_*$ was used to define the splitting of equation
  \eqref{eq:abs_equation}, then the system has at most one invariant
  measure.
\end{theorem}

Note: It is worth mentioning that existence of an invariant measure in
our setting is usually straight forward. For instance, if the set $\{ u :
V(u) \leq M\}$ is  precompact for all $M$ then the result follows
easily by the standard Krylov--Bogoljubov construction of extracting a
convergent subsequence from the empirical measures obtained by
time--averaging. See for instance \cite{b:ChowKhasminskii98} for the SPDE
setting or \cite{b:CornfeldFominSinai82} for general discussions.

We begin the proof of Theorem \ref{thm:abs_ergodic} by proving the
analog of Lemma \ref{l:toy_marginals} from the discussion of the toy
model. In fact, we will only deduce part of it from our existing
assumptions leaving the remainder still as an assumption.

\begin{lemma} \label{l:abs_equiv}Consider the solution to equation
  \eqref{eq:abs_equation} under the assumptions of Theorem
  \ref{l:abs_equiv}. For $\ell_0 \in \Xl$ and $h_0,\tilde h_0 \in \Xh$
  and almost every $\eta$, $Q_{[0,\infty)}^\eta(\ell_0,h_0,\ccdot)$ is
  equivalent to $Q_{[0,\infty)}^\eta(\ell_0,\tilde h_0,\ccdot)$.
\end{lemma}

\bpf[Proof of Lemma \ref{l:abs_equiv}] Again we begin
by essentially localizing to a fixed $B_n$. However, we need to
pick a set of paths in  $C(0,\infty;\Xl)$. Fixing $\ell_0, h_0, \tilde
h_0$ and $\eta$, we define
\begin{multline*}
  B_n'=\Big\{ \ell_{[0,\infty)} : u, \tilde u \in B_n \mbox{ where }
  u(s)=\ell(s) + \Phi_{s}^\eta(\ell_{[0,s]},h_0), \\ \tilde u(s)=\ell(s)
  + \Phi_{s}^\eta(\ell_{[0,s]},\tilde h_0) \Big\}
\end{multline*}
Then for $A\subset C(0,\infty;\Xl$), define
\begin{align}
  \label{eq:Q0tB_n}
  Q_{[0,t]}^\eta(\ell_0,h_0,A;B_n')&=\PP\Big\{ u_{[0,t]} \in A;
  \Pl u_{[0,t]} \in B_n' \Big|
  u(0)=(\ell_0,h_0)\Big\} 
\end{align}
We now compare $Q_{[0,t]}^\eta(\ell_0,h_0,\ccdot;B_n')$ to
$Q_{[0,t]}^\eta(\ell_0,\tilde h_0,\ccdot;B_n')$.  Again we compare the
measure using Lemma \ref{l:compareMeasures} from the appendix. By
restricting to $B_n'$, we ensure that both $u$ and $\tilde u$ stay in
$B_n$. Hence, the first part of Lemma \ref{l:abs_contraction} combined
with the second estimate in Assumption \ref{a:contractive} produces
\begin{multline*}
  \normX{\Pl G(u(t))-\Pl G(\tilde u(t))}^2 \leq c_4\left[ 1 +
    2(C_1t+n[1+\log(1+t)])^{p_1}\right] \times \\
\normX{ \Phi_{0,t}^\eta(\Pl u_{[0,t]};h_0) -
   \Phi_{0,t}^\eta(\Pl u_{[0,t]};\tilde h_0)}^{p_2}
\end{multline*}
where
 \begin{align}\label{eq:contractionEstimateForEnd}
   \normX{ \Phi_{0,t}^\eta(\Pl u_{[0,t]};h_0) - \Phi_{0,t}^\eta(\Pl
     u_{[0,t]};\tilde h_0)}^2 \leq\normX{h_0 - \tilde h_0}^2 e^{nc_*}
   (1+t)^{nc_*} e^{-\frac12\gamma_* t} \ .
  \end{align}
  Defining $\sigma_{min}^2=\min_{|k| < N_*} |\sigma_k|^2$, the
  previous two estimates imply that
\begin{align}\label{eq:diffGirsonovEstimate}
  \int_0^\infty  \frac1{\sigma_{min}^2} \normX{\Pl G(u(t))-\Pl
  G(\tilde u(t))}^2 dt <D_*<  \infty
\end{align}
for some $D_*$ uniformly on $B_n'$. Using Lemma
\ref{l:compareMeasures}, we conclude that
$Q_t^\eta(\ell_0,h_0,\ccdot;B_n')$ is equivalent to
$Q_t^\eta(\ell_0,\tilde h_0,\ccdot;B_n')$. As in the previous part,
since both $u$ and $\tilde u$ are in $\cup B_n$ with probability one,
we conclude that $Q_t^\eta(\ell_0,h_0,\ccdot)$ is equivalent to
$Q_t^\eta(\ell_0,\tilde h_0,\ccdot)$.

Looking back on the above proof, we seen that there was a great deal
of uniformity in the estimates. When comparing
$Q_{[0,t]}^\eta(\ell_0,h_0,\ccdot;B_n')$ to
$Q_{[0,t]}^\eta(\ell_0,\tilde h_0,\ccdot;B_n')$, we see that for all
$(\ell_0,h_0),(\ell,\tilde h_0)$ in a bounded ball, we can choose the
same $D_*$. From Lemma \ref{l:compareMeasures} in the appendix, we get
the following result
\begin{lemma}\label{l:RNuniformEstimate}
  For any $M$, there exists a $D_*$ so that if
  $\normX{\ell_0+h_0},\normX{\ell_0+\tilde h_0}
  \leq M$ then 
  \begin{align*}
 \EE\left[
 \frac{dQ_{[0,t]}^\eta(\ell_0,h_0,\ccdot;B_n')}{dQ_{[0,t]}^\eta(\ell_0,\tilde
 h_0,\ccdot;B_n')}\right]^p \leq D_*^{p(p-1)}
  \end{align*}
  for all $p>0$.
\end{lemma}
\bpf[Conclusion of the Proof of Theorem \ref{thm:abs_ergodic}]
In light of Lemma \ref{l:abs_equiv} and \ref{l:abs_contraction} the
result follows from Theorem \ref{thm:ergodic}.
\epf

We now address Assumption \ref{a:badAssumptions}. In the case of the
stochastic Navier Stokes equations it is implied without further
assumptions by the techniques used to prove Theorem
\ref{thm:densityPDE} from section \ref{sec:hypoelliptic} since all
of the directions in $\Xl$ are directly forced. (Theorem
\ref{thm:densityPDE} does not address the question of the marginal
with respect to $\eta$. However Theorem
\ref{thm:densityPDE} follows from the fact that the Malliavin
covariance matrix restricted to $\Xl$ is almost surely
invertible. This does imply the result for the marginals. See
\cite{b:MattinglyPardoux03Pre}.)

The same techniques should apply to most SPDEs of interest with
additive noise. However, since an abstract version of the techniques
in \cite{b:MattinglyPardoux03Pre} is not written, we refrain from
making any claims. There is however another approach. Though it is
rather adhoc and in the author's opinion and ``not the correct way,'' it
is sufficient in many contexts. The basic idea is to compare the
measures induced on $C(0,t;\Xl)$ by the process of interest and some
well understood process both starting from the same point. This is
done using Girsonov's theorem. Then the time $t$ marginals of the well
understood process starting from two different points are compared. By
stringing the estimates together and making some additional
assumptions needed to control the ``high'' modes in  equation
\eqref{eq:highMode}, one can prove Assumption \ref{a:badAssumptions}. A
particularly simple version of this was done in the toy model of the
previous section. For more complicated versions see
\cite{b:EMattinglySinai00,b:ELui02,b:Mattingly02,b:BakhtinMattingly03Pre}.
\cite{b:BakhtinMattingly03Pre} has a relatively crisp version of the argument.

\section{Exponential Mixing and Coupling}
\label{sec:expMixing}

In this section, we expand the simple uniqueness results, given
earlier in the paper, by giving a rate of convergence. The proof will
be based on a coupling argument and is closer in packaging to the
author's first proof of basic ergodicity which were presented in 
seminar talks\footnote{Stanford and Berkeley probability seminars
  November and December 1999.}. We will measure the rate of convergence
of \eqref{eq:abs_equation} using the following metric. For any two
measures $\mu_1$ and $\mu_2$ on $\X$ define
\begin{align*}
  \|\mu_1 - \mu_2\|_* = \sup_{\phi \in \mathcal{G}_*} \int
  \phi(x) \mu_1(dx) - \int \phi(x) \mu_2(dx)
\end{align*}
where $\mathcal{G}_*$ is the set of all measureable functions $\phi:\X
\rightarrow \RR$ with $ |\phi(u)| \leq 1$ for all $u \in \X$ and
$|\phi(\ell + h)- \phi(\ell + \tilde h)| \leq \normX[big]{h - \tilde
  h}$ for all $\ell \in \Xl$ and $h,\tilde h \in \Xh$. Notice that the
$\|\ccdot\|_*$ norm dominates the Wasserstein or Kantorovich distance
for measures but is weaker than the total variation norm.  In the
definition of $\mathcal{G}_*$, we could have also used test function
which were $Lip_\alpha$ on $\Xh$, with $\alpha >0$, and all of the
theorems below would still hold.

Now we make the following assumption which is a more qualitative
version of Assumption \ref{a:badAssumptions}. It amounts to continuity
in the initial condition of the density induced on $\Xl$ at time $t$.
\begin{assumption} \label{a:bad2}
  Fix any $t >0$. For any $M_0$, there exist a positive  $\delta$ and
  $\Omega' \subset \Ph \Omega$ so $\PP(\eta \in \Omega') > \delta$ and
  for any $\eta \in \Omega'$ and $u_0^{(i)} \in \X$, $i=1,2$, with
  $V(u_0^{(i)}) \leq M_0$ we have $\frac12\|Q_t^\eta(u_0^{(1)},\ccdot) -
  Q_t^\eta(u_0^{(2)},\ccdot)\|_{TV} < 1- \delta$ .
\end{assumption}
See Appendix \ref{sec:couplingEstimatesAppendix} for the definition of
$\|\ccdot\|_{TV}$ which may differ by a factor of 2 from some definitions.
Again this estimate can be obtained in a number of ways. For the SNS
it was obtained by comparing, in a quantitative fashion, the total
variation distance between the time $t$ marginals and well controlled
reference process (either Brownian motion or the SDE on $\Xl$ obtained
from the Galerkin truncation of the SNS). However the author feels that
this is not the optimal fashion to proceed. It would be better to use the flow
property  and the calculations from \cite{b:MattinglyPardoux03Pre} to
verify this estimate. Since the assumption has only been verified in
specific cases, we leave it as an assumption.

Letting $P_t(u_0,A)=\PP\{ u(t) \in A | u(0)=u_0\}$ where $A \subset
\X$, we have the following result whose proof give in the sections
which follow. Stronger results using norms allowing test functions
which grow are also possible by the methods presented here. Corollary
\ref{c:extra} at the end of the section gives a simple, suboptimal
example. See \cite{b:MeTw93} or \cite{b:MattinglyStuartHigham00} for
examples to the type of stronger statements which should be possible.
However \cite{b:MeTw93,b:MattinglyStuartHigham00} does not apply to
our setting.
\begin{theorem}\label{thm:expConvergence} If Assumption
  \ref{a:Lyop}, \ref{a:contractive} and \ref{a:bad2} hold, then there
  exists fixed positive constants $K$ and $\gamma$ so for all
  $u_0^{(i)} \in \X$ (possibly random, but adapted to the filtration
  at time zero) (see notes below.)
  \begin{align*}
 \|P_t(u_0^{(1)},\ccdot) - P_t(u_0^{(1)},\ccdot)\|_*
  \leq K\big[1+ \EE V(u_0^{(1)}) +  \EE V(u_0^{(2)}) \big] e^{-\gamma t}     
  \end{align*}
\end{theorem}
 
We give the proof of this theorem in the next sections.  In
\cite{b:Mattingly02} a general theorem, ensuring exponential mixing in
a wide class of problems,  was given and the conditions
were verified for the SNS. However, given the estimates of the
previous section the exact same analysis applies to equation
\eqref{eq:abs_equation} when $\eta=\Ph W=0$. In that paper, the case
$\eta=\Ph W\not=0$ was discussed in another setting.  A straight
forward modification of the techniques from that paper yields the
extension to $\eta=\Ph W\not=0$.  Kuksin and Shirikyan were the first
to consider exponential mixing for the SNS in the case when
$\eta\not=0$ \cite{b:KuksinShirikyan02a}; however, their norm is
slightly weaker. The norm we give here gives total variation
convergence on a subset of the space which dictates the asymptotic
behavior, namely $\Xl$. This allows on to use standard mixing results
to get law of large numbers, central limit theorems, and other
results. With additional work this also possible directly in the framework of
\cite{b:KuksinShirikyan02a} or \cite{b:Hairer02}. See \cite{b:Shirikyan02}

In \cite{b:Mattingly02}, the case $\eta \neq 0$ was considered in a
simple map example and we see here that those ideas extend to the SPDE
context. In \cite{b:BricmontKupiainenLefevere02}, exponential
convergence was proven but without the explicit dependence on the
initial condition.  That paper along with \cite{b:Mattingly02} were
the first proofs of exponential convergence of the SNS with white in
time forcing. In the kicked case exponential convergence was given in
\cite{b:MasmoudiYoung02,b:KuksinPiatnitskiShirikyan02}. The first of
these also considers the the case where the system is strongly
dissipative as in Theorem \ref{thm:largeNu}. In \cite{b:Hairer02},
exponential convergence for a reaction diffusion equation was proved
by bringing the paths together asymptotically using a coupling
construction inspired by \cite{b:Mattingly02}. (Both
\cite{b:Mattingly02} and \cite{b:BricmontKupiainenLefevere02} were
delayed considerable in the review process, and hence,
\cite{b:Hairer02} appeared first.)
 
To state a slightly stronger result, for any weighing function $R: \X
\rightarrow [0,\infty)$ define $\|\mu_1 - \mu_2\|_{R*} = \sup_{\phi
  \in \mathcal{G}_R*} \int \phi(x) \mu_1(dx) - \int \phi(x) \mu_2(dx)$
where $\mathcal{G}_{R*}$ is the set of all measureable functions
$\phi:\X \rightarrow \RR$ with $ |\phi(u)| \leq R(u)$ for all $u \in
\X$ and $|\phi(\ell + h)- \phi(\ell + \tilde h)| \leq\big[1+ R(\ell+h)
+ R(\ell+\tilde h)\big] \normX{h - \tilde h}$ for all
      $\ell \in \Xl$ and $h,\tilde h \in \Xh$.
\begin{corolary}
  \label{c:extra} In the same setting as Theorem
    \ref{thm:expConvergence}, for any $\frac1q +\frac1p =1$ with $q,p >1$ 
    \begin{multline*}
 \|P_t(u_0^{(1)},\ccdot) - P_t(u_0^{(1)},\ccdot)\|_{R*}
 \leq[1+(\EE R(u(t))^q)^\frac1q+(\EE R(\tilde u(t))^q)^\frac1q]\times\\
K'\big[1+ \EE V(u_0^{(1)}) +  \EE V(u_0^{(2)}) \big]^\frac1p
 e^{-\gamma' t}    \end{multline*}
where $\gamma'$ and $K'$  are positive constants depending on $p$ and $q$.
\end{corolary}
Notice that if $R(x)=V(x)$ then the assumptions of the corollary are
satisfied and $\EE V(u(t))^q \leq K''[1 + \EE V(u_0)^q]$ for some
$K''$ as $V^q$ is also a Lyapunov function. This Lemma is suboptimal
as the right hand side does not scale linearly in $V(u_0)$ so a
convenient operator norm is not induced. See
\cite{b:MattinglyStuartHigham00} for ideas, from the Markov setting, which likely could overcome
this difficnency.
\subsection{Deconstruction and Reconstruction}
\label{sec:deconRecon}

We begin with an overview of the coupling construction. The idea is to
factor the measure induced on $C(0,\infty;\Xl)$ starting from $u_0$
and $\tilde u_0$ and build a process on $C(0,\infty;\X) \times
C(0,\infty;\X)$ so that the marginals are distributed as a process
started from $u_0$ and $\tilde u_0$ respectively and so that $\Pl u(t)
= \Pl \tilde u(t)$ with positive probability. There is the added
complication that we need to also have the processes use the same
realization of $\eta=\Ph W$ and that we need to localize the
trajectories to the nicely growing and averaging paths so that  $\Ph u(t)
- \Ph \tilde u(t)$ will converge to zero at a controlled rate. We begin
with the localization.

The $B_n$ defined in section \ref{sec:fluctuationsGeneral} were
sufficient for localizing to prove uniqueness. They also showed how
typical paths stayed in a logarithmic envelope about the average
behavior. However the probability from deviating from a given $B_n$
after time $t$ decays slowly. To prove exponential convergence, we now
localize with sets from which it becomes exponentially unlikely to
deviate over time. For positive $M$ define
\begin{multline}
  \label{eq:BExp}
   B(M)=\Bigl\{ u\in C(0,\infty; \X):
  V(u(t)) +(1-\epsilon_*)\int_0^t U(u(s))ds -
    V(u_0) \\ \leq M+ C_1(1+\epsilon_*)t \ \mbox{ for all $t\geq 0$ }\Bigr\} 
\end{multline}
The constant $\epsilon_*$ is chosen so that
$c_1-\frac{1+\epsilon_*}{1-\epsilon_*}C_1(c_2+c_3)= \frac12\gamma_*$.
Recall that $\gamma_*=c_1-C_1(c_2+c_3)$ was assumed positive. Clearly 
$\PP\{ u \in \cup_{M=1}^\infty B(M) \} =1$ and \eqref{eq:expMartEstimate}, $\PP\{ u_{[0,t]}
\in \Pi_{[0,t]}B(M); u \not \in B(M) \}$ decays exponential in $t$.
Furthermore given the choice $\epsilon_*$, Lemma
\ref{l:abs_contraction} (part one) holds with $B_n$ replaced with $B(M)$
and different constants on the right hand side of the decay estimate.
Precisely, if $u^{(i)}_{[0,T)} \in \Pi_{[0,T)}B(M)$ ($i=1,2$) where
$u^{(i)}(t)=\ell(t)+ \Phi_{0,t}^\eta(\ell_{[0,t]};h_0^{(i)})$ then for $t
\in (0,T]$
\begin{equation}
  \label{eq:CoupContraction}
    \normX{ \Phi_{0,t}^\eta(\ell_{[0,t]};h_0^{(1)}) -
   \Phi_{0,t}^\eta(\ell_{[0,t]}; h_0^{(2)})}^2 \leq\normX{h_0^{(1)} -
   h_0^{(2)}}^2 e^M e^{-\frac12\gamma_* t} \ .
\end{equation}

Fix some $M$. For every $(u_0^{(1)},u_0^{(2)},\eta) \in
\X\times\X\times\Ph\Omega$, we define
\begin{multline*}
  \bar B_{[0,n]}(u_0^{(1)},u_0^{(2)},\eta)=\Big\{ \ell_{[0,n]} \in
  \Pl\X_{[0,n]} : u_{[1,n]}^{(i)} \in \Pi_{[0,n-1)}B(M),
  \normX{u^{(i)}(1)}^2\leq M \\ \mbox{ where } u^{(i)}(s)=\ell(s) +
  \Phi_{s}^\eta(\ell_{[0,s]},\Ph u_0^{(i)}) \Big\}
\end{multline*}
and define for $A \subset \Pi_\ell\X_{[0,n-1)}$ the measure
\begin{align*}
  Q_{[1,n]}^\eta(u_0, A ;  \bar B_{[0,n]})&= \PP\big\{ \Pl u_{[1,n]}
  \in A  \mbox{ and }   \Pl u_{[0,n]} \in  \bar
  B_{[0,n]}\big| u(0)=u_0,\mathcal{F}_{[0,n]}^\eta \ \big\} \\
  &=\EE\{\one_A( \Pl u_{[1,n]}) \one_{ \bar B_{[0,n]}}( \Pl
  u_{[0,n]}) | u(0)=u_0, \mathcal{F}_{[0,n]}^\eta\} .
\end{align*}
where $ \bar B_{[0,n]}= \bar B_{[0,n]}(u_0^{(1)},u_0^{(2)},\eta)$.
Hence $Q_{[1,n]}^\eta(u_0, A ; \bar B_{[0,n]})$ is the measure of
paths so that $\Pl u_{[0,n]} \in \bar B_{[0,n]}$ and $ \Pl u_{[1,n]}
\in A$ if one conditions to start from $u_0$ and use noise
realizations $W$ so $\Pl W=\eta$. Of course, it is not a probability
measure as it does not have total mass one.

Given any two measures $\mu_1$ and $\mu_2$, one can always write them
as a density against a common third measure. That is
$\mu_i(dx)=f_i(x)\mu_3(dx)$ for $i=1,2$. We define the measures
$\mu_1\wedge \mu_2$ and $\mu_1 - \mu_2$ respectively by the densities
$(f_1(x)\wedge f_2(x))\mu_3(dx)$ and $(f_1(x)-f_2(x))\mu_3(dx)$. It is
easy to see that this definition is independent of the choice of $\mu_3$. If
$\mu_1$ does not dominate $\mu_2$ for all measurable sets then the
second measure is a signed measure. See Appendix
\ref{sec:couplingEstimatesAppendix} more explination and the realtion
to the total variation norm, which we denote by $\|\ccdot\|_{TV}$.

With this notation define
\begin{align*}
  \Psi^\eta_n(u_0^{(1)},u_0^{(2)},\ccdot)=  Q_{[1,n)}^\eta(u_0^{(1)}, \ccdot
  ;  \bar B_{[0,n)})\wedge  Q_{[1,n)}^\eta(u_0^{(2)}, \ccdot ;  \bar B_{[0,n)})
\end{align*}
where again $\bar B_{[0,n)}= \bar
B_{[0,n)}(u_0^{(1)},u_0^{(2)},\eta)$. Next for $A \subseteq
\X_{[0,n)}$ and $\ell \in \Pl \X_{[0,n-1)}$ define
\begin{align*}
  H_n^\eta(u_0,A | \ell )= \PP\{ u_{[0,n)} \in A | u(0)=u_0, \Pl
  u_{[1,n)} = \ell_{[0,n-1)},\mathcal{F}_{[0,n)}^\eta \} .
\end{align*}
In words $H_n^\eta(u_0,\ccdot | \ell )$ is the measure induced on
$\X_{[0,n)}$ by paths $u_{[0,n)}$ conditioned to start at $u_0$, use
noise realization $\eta$, and such that $u(s+1)=\ell(s)$ for $s \in
[0,n-1)$.

Next we define the two families of measures $\r_n$ and $\s_n$,
$n\in\{1,2,\dots,\infty\}$, which will be critical in our
construction. They will both be measures on $\Pl\X_{[0,n)} \times
\Pl\X_{[0,n)} \times \Ph \Omega_{[0,n)}$ with
$n\in\{1,2,\dots,\infty\}$. In general, we will use bold letters to
denote measures on such spaces and capital bold letters for
probability measures on such spaces. Define
\begin{multline*}
  \s_n(u_0,\tilde u_0, du \times d\tilde u \times d\eta) \\
= \int_{\Pl \X_{[0,n-1)}} [
  H_n^\eta(u_0,du | \ell) \times  H_n^\eta(\tilde u_0,d\tilde u | \ell)
  ]\Psi^\eta_n(u_0,\tilde u_0,d\ell) \times \PP(d\eta)
\end{multline*}
 and
 \begin{gather*}
   \r_{n+1}(u_0,\tilde u_0, du \times d\tilde u \times
   d\eta)=\big[\mathbf{P}_1\s_n -\s_{n+1}\big](u_0,\tilde u_0, du \times
   d\tilde u \times d\eta) .
 \end{gather*} Here $ \P_1\s_n$ is the measure on $\Pl\X_{[0,n+1)} \times
\Pl\X_{[0,n+1)} \times \Ph \Omega_{[0,n+1)}$ obtained by first stepping
with $\s_n$ and then with 
\begin{gather*}
  \P_n(u_0,\tilde u_0,du_{[0,n)} \times d \tilde u_{[0,n)} \times
  d\eta_{[0,n)})= P_{[0,n)}^\eta(u_0,d u_{[0,n)}) \times
  P_{[0,n)}^\eta(\tilde u_0 d
  \tilde u_{[0,n)}) \times \PP(d\eta_{[0,n)})
\end{gather*}
where for $A \subset \X_{[0,n)}$ and $P_{[0,n)}^\eta(u_0,A)=\PP\big\{
u_{[0,n)} \in A \big| u(0)=u_0, \mathcal{F}^\eta_{[0,n)}\big\}$. That
is to say,
\begin{multline*}
 \P_1 \s_n (u_0,\tilde u_0, du_{[0,n+1)} \times d\tilde
  u_{[0,n+1)}\times d\eta_{[0,n+1)})= \s_n(u_0,\tilde u_0, du_{[0,n)}
  \times d\tilde u_{[0,n)} \times d\eta_{[0,n)}) \times\\\P_1(u(n),\tilde u(n),
  du_{[n,n+1)} \times du_{[n,n+1)} \times  d\eta_{[n,n+1)}).
\end{multline*}
Define $\r_1(u_0,\tilde u_0,\ccdot) =\P_1(u_0,\tilde
u_0,\ccdot)-\s_1(u_0,\tilde u_0,\ccdot)$.  Since $ \P_1
\s_n(u_0,\tilde u_0,A) \geq \s_{n+1}(u_0,\tilde u_0,A)$ by
construction for all measurable sets $A$, $\r_n$ is a standard measure
and not a signed measure. Lastly, we define
\begin{align*}
  \rho_n(u_0,\tilde u_0)= \s_n(u_0,\tilde u_0, \X_{[0,n)} \times
  \X_{[0,n)} \times \Ph \Omega_{[0,n)} )
\end{align*}
for $n>0$ (including $n=\infty$) and $\rho_0=1$ and the probability
transition kernels
\begin{align*}
  \S_n(u_0,\tilde u_0,\ccdot)=\frac{\s_n(u_0,\tilde
    u_0,\ccdot)}{\rho_n(u_0,\tilde u_0)} \qquad \R_n(u_0,\tilde
  u_0,\ccdot)=\frac{\r_n(u_0,\tilde u_0,\ccdot)}{\rho_{n-1}(u_0,\tilde
    u_0)-\rho_{n}(u_0,\tilde u_0)} .
\end{align*}
If the denominator is zero in either of the above definitions, we set
the corresponding measure to the zero measure.  Observe that
 \begin{align*}
   \rho_n(u_0,\tilde u_0)=& \EE \Psi^\eta_n(u_0,\tilde u_0,\Pi_\ell \X_{[0,n-1)}) \\
   =&1-\frac12 \EE \| Q_{[1,n)}^\eta(u_0, \ccdot
   ;  \bar B_{[0,n)}) - Q^\eta_{[1,n)}(\tilde u_0, \ccdot ;  \bar B_{[0,n)})
   \|_{TV} > 0, 
 \end{align*}
 where $ \bar B_{[0,n)}= \bar B_{[0,n)}(u_0,\tilde u_0,\eta)$. This
holds even for $n=\infty$, since for all $n$ the measures are
absolutely continuous for almost every $\eta$. This can be seen by the
same calculations as in the proof of Lemma \ref{l:abs_equiv} coupled
will Lemma \ref{l:Equiv}.  Also observe that $\rho_n(u_0,\tilde u_0)
\geq \rho_{n+1}(u_0,\tilde u_0)$. Thus, we have
\begin{align*}
  1=\rho_0 \geq \rho_1 \geq \cdots \geq \rho_\infty 
\end{align*}
For all $M$ sufficiently large, we will see in Lemma \ref{l:rhos} that
$\rho_\infty(u_0,\tilde u_0) > 0$ for all $u_0, \tilde u_0$ with
$V(u_0), V(\tilde u_0) \leq M_0$. 

From the properties of $\s_n$ and $\r_n$, one had $\P_1(u_0,\tilde
u_0,\ccdot) = \s_1(u_0,\tilde u_0,\ccdot) + \r_1(u_0,\tilde
u_0,\ccdot)$ and
 \begin{align*}
   \P_2 &\eqM \P_1\s_1 +  \P_1\r_1 \eqM \s_2 + [\P_1\s_1 - \s_2] + \r_1\\
   &\eqM \s_2 + \r_2 + \r_1
 \end{align*}
 where we have suppressed the dependence of the kernels on the initial
 conditions $u_0$ and $\tilde u_0$.  By $\eqM$ we mean that the two
 measures have the same relevant marginals. More precisely if we
 consider the kernel at the point $(u_0,\tilde u_0)$, the joint
 distribution of the first and last coordinate of both sides is
 $P_{[0,2)}^\eta(u_0,\ccdot)\times \PP(d\eta)$ and the joint
 distribution of the second and last coordinate of both sides is
 $P_{[0,2)}^\eta(\tilde u_0,\ccdot)\times \PP(d\eta)$.  Continuing
 along this line and normalizing the measures to probability measures,
 produces the following version of the factoring lemma from
 \cite{b:Mattingly02}.
 \begin{align}
   \label{eq:factorization}
   \P_\infty(u_0,\tilde u_0,\ccdot) \eqM \rho_\infty \S_\infty +
   \sum_{n=1}^\infty [\rho_{n-1} -\rho_n] \P_\infty\R_n
 \end{align}
 where $\P_\infty \R_n$ is analogous to $ \P_1 \s_n$ from above. On
 the right hands side we have suppressed the dependence on $u_0$ and
 $\tilde u_0$ in the interest of space. That is $ \S_\infty=
 \S_\infty(u_0,\tilde u_0,\ccdot)$, $\rho_{n-1}=
 \rho_{n-1}(u_0,\tilde u_0)$ and so forth.
 
 Such a factorization of the futures was also fundamental to the
 results in \cite{b:Hairer02}.  Since the project of this measure onto
 the first and last coordinate of both sides equals
 $P^\eta_{[0,\infty)}(u_0,\ccdot) \times \PP(d\eta_{[0,\infty)})$ and
 the projection on the second and last coordinate of both sides is
 $P^\eta_{[0,\infty)}(\tilde u_0,\ccdot) \times
 \PP(d\eta_{[0,\infty)})$, we have built a representation of two
 copies of the process which both use the same $\eta$. The first is
 distributed as a solution starting from $u_0$ and the second as a
 solution starting from $\tilde u_0$.  This representation has the the
 following importantly feature. There exists a set $A \subset
 \X_{[0,\infty)} \times \X_{[0,\infty)} \times \Ph
 \Omega_{[0,\infty)}$ so that $\S_\infty(A)=1$ and if $(u,\tilde
 u,\eta) \in A$ then $\Pl u(s) =\Pl \tilde u(s)$ for all $s \geq 1$,
 $u_{[1,\infty)}, \tilde u_{[1,\infty)} \in B(M)$, and $u, \tilde u$
 are solutions for some noise realizations $W$ and $\tilde W$ so $\Ph
 W= \Ph \tilde W$. These are precisely the conditions needed to apply
 the contractive estimates from section \ref{sec;genContraction}.

 This factorization states that drawing from $ \P_\infty$ is
 equivalent, as far as either $u(t)$ or $\tilde u(t)$, is concerned,
 to drawing from $\S_\infty$ with probability $\rho_\infty$ and
 $\P_\infty\R_n$ with probability $\rho_{n-1} -\rho_n$. Of course, we
 have built in useful correlations between the two processes.  Also
 notice that $\P_\infty$ appears on the left hand side, so the
 factorization can be iterated.

\subsection{Estimates on the $\rho$'s}
\label{sec:rhos}

The following estimates on the $\rho$ are the principle information
needed to prove the exponential mixing, other that the Lyapunov
structure which will be described in the next section. The first
estimate is enough to imply mixing. The fact that the spacing between
the $\rho$'s decays exponentially, combined with the exponential tails
of the return time to the set $\bC$ defined in the following section,
give the exponential mixing rate.

\begin{lemma}\label{l:rhos} In the setting of Theorem \ref{thm:expConvergence},
  let $B(M)$ be the set used to define $\bar B_{[0,n)}$ in the previous
  section. For any $M_0> 0$ the following estimates hold for all $M$
  large enough:
  \begin{enumerate}
  \item There exists a positive constant $\rho_\infty^*$, depending on
  $M$ and $M_0$, so that
  \begin{align*}
    \inf_{u_0^{(i)}:V(u_0^{(i)}) \leq M_0} \rho_\infty(u_0^{(1)},u_0^{(2)}) \geq
    \rho_\infty^* >0
  \end{align*}
\item There exist positive constants $K_1$ and $\gamma_1$, also
  depending on $M$ and $M_0$, so that for all $u_0^{(i)}$, $i=1,2$
  with $V(u_0^{(i)}) \leq M_0$, 
  \begin{align*}
  \rho_{n}(u_0^{(1)},u_0^{(2)}) -\rho_{n+1}(u_0^{(1)},u_0^{(2)}) \leq K_1
  \exp(-\gamma_1 n) .  
  \end{align*}
\end{enumerate}
\end{lemma}
The proof of this lemma will be given in section \ref{sec:rhosProof} .
\subsection{Consequences of the Lyapunov Structure}
\label{sec:ConCLyop}

 We now make a modification in the presentation relative to
 \cite{b:Mattingly02} which is greater than notational (but still
 mainly cosmetic\label{er:cosmetic}). We want to iterate the expansion
 \eqref{eq:factorization}. However we will only have nice control over
 the $\rho_n$'s for $u_0, \tilde u_0$ in a particular subset of the
 phase space. Hence, we modify the expansion to include the steps
 needed to return to this subset.

As already mentioned under Assumption \ref{a:Lyop}, a lemma
analogous to Lemma \ref{l:uniformEnergyEnstrophy} holds for the
Lyapunov function $V$. From this it is straight forward that there
exists an $\alpha \in(0,1)$ so that $\EE\{ V(u(t+1)) |
\mathcal{F}_{t}\} \leq \alpha V(u(t)) +C_1$. Hence, if we define
$\V(u,\tilde u)= V(u)+V(\tilde u)$ then 
\begin{align*}
  \EE\{ \V(u(t+1),\tilde u(t+1)) | \mathcal{F}_{t}\} \leq \alpha
  \V(u(t),\tilde u(t)) +2C_1 .
\end{align*}
We define the set $\bC=\{ (u,\tilde u) : \V(u,\tilde u) \leq \frac{4
  C_1}{\alpha} \}$ and the stopping time
\begin{equation*}
 \tau_C= \inf\{ s \geq 0: s \in \N;  (u(s),\tilde u(s)) \in \bC \} . 
\end{equation*}
Lastly set $M_0=\frac{4 C_1}{\alpha}$ and fix $M$, from the previous
two sections so the conclusions of Lemma \ref{l:rhos} hold. The
importance of this choice of $M_0$ and hence the definition of $\bC$
are given by the following result.
\begin{lemma}\label{l:returnTimes}
   Under Assumption \ref{a:Lyop}, $\PP\{ \tau_{\bC}(u_0,\tilde u_0) > n\}
   \leq K_0\gamma_0^n \V(u_0,\tilde u_0)$ for any $\gamma_0 \in (\alpha,1)$
   and  some positive $K_0=K_0(\gamma)$.
 \end{lemma}
 \bpf[Proof of Lemma \ref{l:returnTimes}] This result can be found
 many places. See for instance Lemma 11.3.9 of \cite{b:MeTw93}, Lemma
 9.3 of \cite{b:MattinglyStuartHigham00} or in the continuous time
 setting and in the context of the SNS Lemma 3.2 \cite{b:EMattingly00}.  
\epf
\subsection{Coupling: A New Representation of the Process}
\label{sec:coupling}
We will define a new presentation of the chain using the factorization
\eqref{eq:factorization}. First however, we modify the factorization
slightly. In light of the previous section, the process $(u(t),\tilde
u(t))$ returns to the set $\bC$ infinitely often at integer times
almost surely. Let $\P_*(u_0,\tilde u_0, \ccdot)$ be the distribution
of $(u_{[0,\tau_C]}, \tilde u_{[0,\tau_C]},\eta_{[0,\tau_C]})$ where
$(u(0),\tilde u(0))=(u_0,\tilde u_0)$. Then $\P_*(u(0),\tilde
u(0),\ccdot)$ is a probability measure on $\chi$ where
\begin{equation*}
  \chi=\bigcup_{k= 0}^\infty \X_{[0,k]} \times  \X_{[0,k]} \times
  \Ph\Omega_{[0,k]} \ .
\end{equation*}
The case $k=0$ is added to cover the situation when $(u_0, \tilde u_0)
\in \bC$ already.  Since we only want to use the previous
factorization for $(u_0,\tilde u_0) \in \bC$, we redefine
$\rho_n(u_0,\tilde u_0)=0$ for $(u_0,\tilde u_0) \not\in \bC$ and set
$\S_n$ equal to the null measure for $(u_0,\tilde u_0) \not\in \bC$.
Hence for $(u_0,\tilde u_0) \not\in \bC$, $\r_1=\P_1$ and all other
$\r_n$ are then the null measure. The result is that for $(u_0,\tilde
u_0) \neq \bC$, the chain takes a step of length one with $u$ and
$\tilde u$ stepping independently.

Returning to the general case $(u_0,\tilde u_0) \in \X \times \X$.
Defining $\R_{n*}=\P_*\R_n$, the factorization
\eqref{eq:factorization} can be rewritten
\begin{equation}
  \label{eq:factorization2nd}
  \P_\infty(u_0,\tilde u_0,\ccdot) \eqM \rho_\infty \S_\infty +
 [1- \rho_\infty]\sum_{n=1}^\infty \frac{[\rho_{n-1} -\rho_n]}{1-
 \rho_\infty } \P_\infty\R_{n*} .
\end{equation}
Again we have suppressed the dependence of the right hand side on the
initial conditions.
Defining
\begin{equation}
  \label{eq:Rinf}
  \R_{\infty *}(u_0,\tilde u_0,\ccdot)= \sum_{n=1}^\infty
 \frac{[\rho_{n-1}(u_0,\tilde u_0) -\rho_n(u_0,\tilde u_0)]}{1- 
 \rho_\infty(u_0,\tilde u_0) } \R_{n*}(u_0,\tilde u_0,\ccdot),
\end{equation}
we consider the chain $X_n=(x_n,\tilde x_n,\eta_n)$ on the state space
$\bar \chi = \chi \cup \big( \X_{[0,k\infty)} \times \X_{[0,\infty)} \times
\Ph\Omega_{[0,\infty)}\big)$ given by taking steps from probability
transition kernel 
\begin{align}
  \label{eq:transition}
  \rho_\infty(u_0,\tilde u_0)\S_\infty(u_0,\tilde u_0,\ccdot) + [1-
  \rho_\infty(u_0,\tilde u_0)] \R_{\infty *}(u_0,\tilde u_0,\ccdot)\ .
\end{align}
We define   
\begin{equation*}
 t_n=\sum_{k=1}^{n-1} |x_k| =\sum_{k=1}^{n-1} |\tilde x_k| 
\end{equation*}
where $|x_k|$ is the length of the trajectory segment $x_k$. $t_n$ is
the time passed in the physical PDE setting after $n$ steps of the
chain have passed. Since the chain adds segments of random length on
each step, $t_n$ is a random quantity.  Similarly associated to $X_n$ is a
trajectory $(u(t),\tilde u(t))$ of the SPDE. It is defined by
$(u(t),\tilde u(t))=\big(x_n(t-t_n), \tilde x_n(t-t_n) \big)$ where
$t_n$ is the unique $t_k$ such that $t_k \leq t \leq t_k+|x_k|$. We
will use both notations depending on which is the most convenient. We
are, of course, only interested in $X_n$ through the step when
$|x_n|=\infty$. This happens the first time a segment is drawn from
$\S_\infty$. For reasons that will be clear, if they are not already,
we refer to this as the ``coupling time.'' We define the stopping time
\begin{align}
  \label{eq:couplingTime}
  \tau=\inf\{ n : |x_n|=\infty\} .  
\end{align}

We pause for a second to notice some of the properties of the chain we
have built.  On the first step if $(u(0),\tilde u(0)) \not\in \bC$, it
takes one step, adding a piece of trajectory of variable, integer
length according to $\P_*(u(0),\tilde u(0),\ccdot)$. Hence, at the end of
this step, the system is in $\bC$. Henceforward each step starts and
ends in $\bC$. With probability $1-\rho_\infty$ the chain draws from
$\R_{\infty *}$. Each of these paths is of finite length. Their
statistics are discussed below. With probability $\rho_\infty$ a path
of infinite length is drawn from $\S_\infty$. After one unit of time,
paths draw from $\S_\infty$ are, by construction, contained in $B(M)$.
In addition by construction, they have norm at time one less than $M$,
use the same $\eta$ increments, and agree on $\Xl$ for $t \geq 1$.
Since at time one the norm is less than $M$, we have an a priori bound
to the separation in the high modes. Thus, if $(v,\tilde v, \eta)$ is
drawn according to $\S_\infty$, then from \eqref{eq:CoupContraction},
$\normX{v(t)-\tilde v(t)}=\normX{\Ph v(t)-\Ph \tilde v(t)} \leq
Me^M\exp(-\frac12\gamma_* t)$.

\subsection{The Heart of the Convergence Result }
\label{sec:convergence result}

We now show how the previous two sections quickly give the needed
estimates to prove Theorem \ref{thm:expConvergence}. For $\phi \in
\mathcal{G}_*$ one has
\begin{align}
  \EE \phi(u(t)) - \phi(\tilde u(t)) &= \EE [\phi(u(t)) -
  \phi(\tilde u(t))][\one_{t_\tau > \frac{t}2 }+\one_{t_\tau \leq \frac{t}2} ] \notag  \\
  & \leq 2 \PP\{ t_\tau > \frac{t}2 \} + Me^M \exp(-\frac14\gamma_* t) .
\end{align}
The first term in the estimate follows from $\phi(u(t)) - \phi(\tilde
u(t)) < 2$ and $\EE \one_{t_\tau > t/2} = \PP\{ t_\tau > t/2 \}$. The
second term follows because for $t> 2t_\tau$ the system has been
following a trajectory drawn from $\S_\infty$ for at least $t/2$ units
of time. Hence, 
\begin{gather*}
        \normX{\Ph u(t) -\Ph \tilde u(t) } \leq
        Me^M\exp(-\frac12\gamma_* \frac{t}{2})
\end{gather*}
as noted in the previous
paragraph. Next observe that $t_\tau(u_0,\tilde u_0) =
\tau_C(u_0,\tilde u_0) + t_\tau(u_{\tau_C},\tilde u_{\tau_C})$ where
$t_\tau(u_0,\tilde u_0)$ and $t_\tau(u_{\tau_C},\tilde u_{\tau_C})$
means the stopping time starting from initial conditions $(u_0,\tilde
u_0)$ and $(u_{\tau_C},\tilde u_{\tau_C})$ respectively. Hence,
 \begin{equation}
\label{eq:enterToCReduction}
   \PP( t_\tau(u_0,\tilde u_0) > n)  \leq \PP( \tau_C(u_0,\tilde u_0)
   > \frac{n}2) + \sup_{(u_0',\tilde u_0') \in \bC} \PP(
   t_\tau(u_0',\tilde u_0') > \frac{n}2) \ .
 \end{equation}
 We know from Lemma \ref{l:returnTimes} that $ \PP( \tau_C(u_0,\tilde
 u_0) > \frac{n}2) $ is exponentially decaying in $n$ with a constant
 which scales linearly with $\V(u_0,\tilde u_0)$. Hence, Theorem
 \ref{thm:expConvergence} would be proven. If we show that $\sup \PP(
 t_\tau(u_0',\tilde u_0') > \frac{n}2)$ decays exponentially in $n$.
 This is done in the next section.

The proof of Corollary \ref{c:extra} follows from  similar reasoning.
\begin{align*}
  \EE \phi(u(t)) - \phi(\tilde u(t)) &= \EE [\phi(u(t)) -
  \phi(\tilde u(t))][\one_{t_\tau > \frac{t}2 }+\one_{t_\tau \leq \frac{t}2} ] \notag  \\
  &\leq \EE R(u(t))\one_{t_\tau > \frac{t}2 } + \EE R(\tilde
  u(t))\one_{t_\tau > \frac{t}2 } \\ & \qquad +  Me^M
  \exp(-\frac14\gamma_* t)[1 + \EE R(u(t)) +\EE R(\tilde u(t))]  \\
&\leq\big[(\EE(R(u(t))^q))^\frac1q+ (\EE(R(\tilde u(t))^q))^\frac1q
  \big]\big[\PP\{ t_\tau > \frac{t}2 \}\big ]^\frac1p \\ &\qquad + Me^M \exp(-\frac14\gamma_* t)[1 + (\EE R(u(t))^q)^\frac1q +(\EE R(\tilde u(t))^q)^\frac1q]  .
\end{align*}
Hence and exponential bound on $\PP\{ t_\tau > \frac{t}2 \}$ will also
complete the proof of the Corollary.
\subsection{Moments of the Coupling Time}
\label{sec:momentsOFCoupling}

We now complete the proof of Theorem \ref{thm:expConvergence} by
providing exponential control of the moments of $t_\tau$. The missing
pieces are the following lemma, which we will proven at the end of this
section and some estimates on the $\rho$'s given in the next section.
\begin{lemma}
  \label{l:momentsRifnity}
  There exist positive constants $\gamma_2$ and $K_2$ so that for all
  $(u_0,\tilde u_0) \in \bC$,  $\EE \exp(\gamma_2 |\R_{\infty
    *}(u_0,\tilde u_0)|) \leq \exp(K_2)$. Where $|\R_{\infty
    *}(u_0,\tilde u_0)|$ is the random variable distributed as the
  length of a segment drawn from $\R_{\infty *}(u_0,\tilde
  u_0,\ccdot)$.
\end{lemma}

Using this lemma we quickly finish the proof of Theorem
\ref{thm:expConvergence}. For any $a \in (0,1)$ and $(u_0,\tilde u_0)
\in \bC$
\begin{align}
  \label{eq:momentsOfCoupling}
  \PP\{ t_\tau(u_0,\tilde u_0) > n\} =& \PP\{ t_\tau > n ; \tau > an
  \} + \PP\{ t_\tau
  > n ; \tau \leq an \} \notag\\
  \leq & \PP\{\tau > an \} + \PP\{ t_\tau
  > n ; \tau \leq an \}\notag\\
  \leq & (1 - \rho_\infty^*)^{\lfloor an \rfloor} + e^{-(\gamma_2 -K_2a)n}
\end{align}
where $\gamma_2$ and $K_2$ are the constants from Lemma
\ref{l:momentsRifnity} and $\rho_\infty^*$ from Lemma \ref{l:rhos}.
The first estimate follows because on each step of the chain there is
at least a $\rho_\infty^*$ chance of drawing from $\S_\infty$.
Accepting the second estimate for a moment, choosing any $a \in (0, 1
\wedge \frac{\gamma_2}{K_2})$ gives exponential decay and completes
the proof. 

To see the second estimate, observe that from  Lemma
\ref{l:momentsRifnity} and the fact that $(u_0,\tilde u_0)
\in \bC$, $\EE \exp(\gamma_2 \sum_{k=1}^{an} |x_k|) \leq \exp( an
K_2)$. Hence one has 
\begin{align*}
  \PP\{ t_\tau > n ; \tau \leq an \} \leq  \PP\{  \sum_{k=1}^{an}
  |x_k| > n ;  \} \leq  e^{-(\gamma_2 -K_2a)n}
\end{align*}

\bpf{Proof of Lemma \ref{l:momentsRifnity}}
Let $|\R_{n*}|$ be the random variable distributed as the length of a
trajectory drawn from $\R_{n*}(u_0,\tilde u_0,\ccdot)$. In what
follows, we suppress the dependence on the initial conditions $(u_0,
\tilde u_0)$ of the $\rho$'s and the transition kernels as we always
consider the same initial conditions.

Define the random variable $\zeta$ as follows by
\begin{equation*}
  \zeta= k \text{ with probability for $k\in \{1,2,\dots\}$}
  \frac{\rho_{k-1}- \rho_{k}}{1-\rho_\infty} .
\end{equation*}
Then $|\R_{\infty *}|$ is distributed as $|\R_{\zeta *}|$. Hence, we
have
\begin{align*}
  \PP\{ |\R_{\infty *}| > n \} &= \PP\{ |\R_{\zeta *}| > n ; \zeta
  >\frac{n}2
  \} + \PP\{ |\R_{\zeta *}| > n ; \zeta \leq \frac{n}2 \}\\
  &\leq \PP\{ \zeta > \frac{n}2 \} + \PP\{ |\R_{\zeta *}| > n ; \zeta
  \leq \frac{n}2 \} .
\end{align*}
The first term decays exponentially by the second part of Lemma
\ref{l:rhos}. This leaves only the last term.
\begin{align*}
  \PP\{ |\R_{\zeta *}| > n ; \zeta \leq \frac{n}2 \} &\leq
  \sum_{k=1}^\frac{n}2 \PP \{ |\R_{k *}| > n \} \frac{\rho_{k-1}-
    \rho_{k}}{1-\rho_\infty}
\end{align*}
Notice that $|\R_{k *}|$ is $k$ plus the time to return to $\bC$
starting from $(u(k),\tilde u(k))$.
Using Lemma \ref{l:returnTimes} and that by definition $[\rho_{k-1}-
\rho_{k}]\R_k=\r_k$ produces
\begin{multline*}
  \PP\{ |\R_{\zeta *}| > n ; \zeta \leq \frac{n}2 \}\\
\leq
  \frac1{1-\rho_\infty} \sum_{k=1}^\frac{n}2  \int \PP \left\{
    \tau_C(u(k),\tilde u(k))  > n - k \right\}\R_k(du_,d\tilde u)
  [\rho_{k-1}- \rho_{k}] \\
 \leq \frac{K}{1-\rho_\infty} \sum_{k=1}^\frac{n}2 \gamma^{n-k}\int
  \V(u(k),\tilde u(k)) \r_k(u_0,\tilde u_0,du_{[0,k)},d\tilde u_{[0,k)}) .
\end{multline*}
By the definition of $\r_k$ one sees that for any measurable set $A$,
$\P_k(u_0,\tilde u_0,A) \geq \r_k(u_0,\tilde u_0,A)$. Since $\V$ is
positive, we have
 \begin{align*}
   \int \V(u(k),\tilde u(k)) \r_k(du_{[0,k)},d\tilde u_{[0,k)}) &\leq
   \int \V(u(k),\tilde u(k))  \P_k(du_{[0,k)},d\tilde u_{[0,k)})\\
   &= \EE\big\{\V(u(k),\tilde u(k))\big| (u(0),\tilde u(0))=  (u_0,\tilde u_0)
   \big\}\\ &\leq K'' \text{ since  $(u_0,\tilde u_0) \in \bC$.}
 \end{align*}
 The uniform bound on the integral used to obtain the last estimate
 comes from a lemma controlling $\V$ completely analogous to Lemma
 \ref{l:uniformEnergyEnstrophy} about the energy of the SNS. It can
 be found in many places. It is simply integrating up the Lyapunov
 estimate in time. See for instance Lemma 9.3 of
 \cite{b:MattinglyStuartHigham00} or Lemma 11.3.9 of
 \cite{b:MeTw93}. Continuing, one has
 \begin{align*}
   \PP\{ |\R_{\zeta *}| > n ; \zeta \leq \frac{n}2 \} & \leq
   \frac{K'}{1-\rho_\infty} \sum_{k=1}^\frac{n}2 \gamma_0^{n-k}\leq
   \frac{K'}{1-\rho_\infty} \gamma_0^{\frac{n}2}
 \end{align*}

 \epf

\subsection{Proof of Lemma \ref{l:rhos}}
\label{sec:rhosProof}

\bpf[Proof of Lemma \ref{l:rhos}] The details of a similar argument
are on page 452 of \cite{b:Mattingly02}.  We begin with the first
statement. For any $M>0$ and $A \subset \Xl$, define
\begin{align*}
  Q^\eta_t(u_0,A ; M)&=\PP( \Pl u(t) \in A ; V(u(t))\leq M | u(0)=u_0,
  \mathcal{F}_{[0,t]}^\eta)\\ &= \EE\{\one_A(\Pl u(t)) \one_{V(u(t))\leq
    M} | u(0)=u_0,\mathcal{F}_{[0,t]}^\eta\}.
\end{align*}
Since $\sup_{u_0: V(u_0) \leq M_0} \EE\{ V(u(t)) \} < \infty$, for all
$M$ sufficiently large one has
\begin{align*}
 \inf_{u_0: V(u_0) \leq M_0} \PP\{
V(u(t)) < M^\frac12 \}> 1 - \delta/10.  
\end{align*}
Hence there exist a $\Omega'' \subset \Omega'$ so $\PP(\eta \in
\Omega'')> \delta/2$ and for all $\eta \in \Omega''$ and $u_0^{(i)}
\in \X$, $i=1,2$, with $V(u_0^{(i)}) \leq M_0$, one has
$\|Q_t^\eta(u_0^{(1)},\ccdot;M^\frac12) -
Q_t^\eta(u_0^{(2)},\ccdot;M^\frac12)\|_{TV} < 1- \delta/2$.

Now define $\Gamma^\eta(u_0,\tilde u_0,\ccdot)=
Q_{[0,\infty)}^\eta(u_0,\ccdot ; u \in B(M))\wedge
Q_{[0,\infty)}^\eta(\tilde u_0,\ccdot ; u \in B(M))$. Then 
\begin{align*}
  \rho_\infty(u_0,\tilde u_0) \geq \frac{\delta^2}4
  \inf_{\substack{u_0,\tilde u_0 \in \X, \Pl u_0=\Pl \tilde
      u_0\\V(u_0),V(\tilde u_0) \leq M^\frac12}}
  \EE\Gamma^\eta(u_0,\tilde u_0,\Pl\X_{[0,\infty)}\times
  \Pl\X_{[0,\infty)}\times\Pl \Omega_{[0,\infty)})
\end{align*}
  Since
\eqref{eq:CoupContraction} holds in this setting, the exact same
calculations as in the proof of the second half of Lemma
\ref{l:abs_equiv} hold producing an estimate identical to Lemma
\ref{l:RNuniformEstimate} with $B_n'$ replaced by $\{u \in B(M)\}$ and
valid for
$u_0, \tilde u_0$ with $V(u_0), V(\tilde u_0) \leq M^\frac12$.
Combining this estimate with Lemma \ref{l:abstractCouplingBound}, we
obtain for any $p>1$
 \begin{align*}
   \inf_{\substack{u_0,\tilde u_0 \in \X, \Pl u_0=\Pl \tilde u_0\\
       V(u_0),V(\tilde u_0) \leq M^\frac12}} \EE\Gamma^\eta(u_0,\tilde
   u_0,\Pl\X_{[0,\infty)}\times \Pl\X_{[0,\infty)}\times\Pl
   \Omega_{[0,\infty)}) \geq \left[ 1-\frac1p\right]
   \frac{C(M)^{\frac{p}{p-1}} }{p^\frac1{p-1} D_*^p} .
 \end{align*}
 where $C(M)=\inf_{u_0:V(u_0)\leq M^\frac12}\EE
 Q_{[0,\infty)}^\eta(u_0, \Pl\X_{[0,\infty)}\times
 \Pl\X_{[0,\infty)}\times\Pl \Omega_{[0,\infty)}; u \in B(M))$ and
 $D_*$ is the constant defined analogously to
 \eqref{eq:diffGirsonovEstimate}. 
 
 Notice that 
 \begin{gather*}
         \EE Q_{[0,\infty)}^\eta(u_0, \Pl\X_{[0,\infty)}\times
        \Pl\X_{[0,\infty)}; u \in B(M)) =\PP\{ u_{[0,\infty)} \in B(M) |
        u(0)=u_0\}.
 \end{gather*}
 Hence for $M$ sufficiently large, for all $u_0$ with
 $V(u_0) \leq M^\frac12$ there exists a set $\Omega_0''' \subset
 \Pl\Omega$ so that $\PP(\eta \in \Omega_0''') > 1- \delta/100$ and
 for all $\eta \in \Omega_0'''$ $P_{[0,\infty)}^\eta( u_0, B(M)) >
 1-\delta/100$. Hence $C(M) \geq (1-\delta/100)^2$. This completes the
 first claim.

 Now consider the second claim. Setting $\Y_n=\Pl\X_{[0,n)}\times
 \Pl\X_{[0,n)}\times\Pl \Omega_{[0,n)}$, notice that
 $\rho_{n-1}(u_0,\tilde u_0)-\rho_n(u_0,\tilde u_0)=\r_{n}(u_0,\tilde
 u_0,\Y_n)=[\P_1\s_{n-1} - \s_n](u_0,\tilde u_0,\Y_n)$. From this we
 see that $\rho_{n-1}-\rho_n$ is the probability of drawing from
 $\s_{n-1}$ but not from $\s_n$. There are two ways this can happen.
 First the trajectory can leave the set $B(M)$ between time $n-1$ and
 $n$. This probability is exponentially small in $n$ by the
 construction of $B(M)$ and the estimate in \eqref{eq:BExp}. The
 second way is to draw from the part of distribution contained in
 $B(M)$ between time $n-1$ and $n$ but not in the common part of the
 two $Q^\eta$ distributions.  Over $[0,n-1]$ trajectories $(u,\tilde u,\eta)$
 are drawn from $\s_{n-1}$. Hence almost every trajectory has the
 properties that $\Pl u_{[1,n-1]} = \Pl \tilde u_{[1,n-1]}$ and both
 are in $\Pi_{[0,n-1]}B(M)$. The contractive property derived
 analogously to \eqref{eq:contractionEstimateForEnd} then implies,
 $\normX{u(n-1)- \tilde u(n-1)} \leq Me^M e^{-\frac12 \gamma_*
   (n-1)}$. Let $\bar B_{[n-1,n]}(u_{[0,n-1)},\tilde
 u_{[0,n-1)},\eta_{[0,n-1)})$ be the paths in $\Pl \X_{[n-1,n]}\times
 \Pl \X_{[n-1,n]} \times \Pl \Omega_{[n-1,n]}$ so that when added to
 $(u_{[0,n-1]},\tilde u_{[0,n-1]},\eta_{[0,n-1]})$ the resulting path
 $(u_{[0,n]},\tilde u_{[0,n]},\eta_{[0,n]})$ is such that $u_{[0,n]},
 \tilde u_{[0,n]} \in \Pi_{[0,n]}B(M)$. (As before the part of the
 trajectory in $\Ph\X_{[n-1,n]}$ has to be reconstructed with the aid
 of $\Phi$.)  Hence $Q^\eta_{[0,1]}(u(n-1),\ccdot; \bar B_{[n-1,n]})$
 and $Q^\eta_{[0,1]}(\tilde u(n-1),\ccdot; \bar B_{[n-1,n]})$ where
 $\bar B_{[n-1,n]}=\bar B_{[n-1,n]}(u_{[0,n-1]},\tilde
 u_{[0,n-1]},\eta_{[0,n-1]})$, are the two distributions which will be
 used to draw the next unit length step. Thus the term we need to
 control is
 \begin{multline*}
   \frac12\EE\| Q^\eta_{[0,1]}(u(n-1),\ccdot; \bar B_{[n-1,n]})-
   Q^\eta_{[0,1]}(\tilde u(n-1),\ccdot; \bar B_{[n-1,n]})\|_{TV} 
   \\ \leq \EE\left(\EE\left\{\left[ \frac{ dQ^\eta_{[0,1]}(u(n-1),\ccdot)}{
       dQ^\eta_{[0,1]}(\tilde u(n-1),\ccdot)} -1 \right]^2\one_{\bar
       B_{[n-1,n]}} \Bigg| \mathcal{F}^\eta \right\} \right)^\frac12\\
 \leq \left(\exp(Ke^{-\frac12\gamma_* n}) - 1 \right)^\frac12
 \end{multline*}
The main estimate comes from the last estimate of Lemma
\ref{l:compareMeasures} applied on the measure conditioned on a fixed
$\eta$ path. The estimate $\exp(Ke^{-\frac12\gamma_* n})$ is the estimate on
the constant $D_*$ used in Lemma \ref{l:compareMeasures}. This
estimate is a consequence of the contractive property noticed above use
to estimate the difference term 
\begin{equation*}
\exp\left(\int_{n-1}^n
\frac1{\sigma_{min}^2} \normX{\Pl G(u(t))-\Pl
  G(\tilde u(t))}^2 dt  \right)
\end{equation*}
 in a fashion analogous to \eqref{eq:diffGirsonovEstimate}.

\section{Other Examples}
\label{sec:otherExamp}
The general assumptions used in the previous example are general
enough to cover a number of SPDEs of interest. A natural second
example where all of our analysis applies is the stochastically forced
Cahn-Allen/Ginsburg-Landau equation
\begin{align}
  \label{eq:GL}
  du(x,t)=\left[\nu \Delta u + u - u^3\right] dt +dW(x,t)
\end{align}
where $W(x,t)=\sum_{\mathcal{K}} e_k(x) \sigma_k \beta(t)$, $\beta_k$
are independent standard Brownian motions, $\sigma_k$ are positive
constants and $e_k$ are the elements of the real Fourier basis $$\{ 1,
\sin(2\pi x), \cos(2\pi x), \sin(4\pi x) ,\cos(4\pi x), \cdots \}.$$
See \cite{b:BakhtinMattingly03Pre,b:ELui02} for the verification of
the assumptions. (Note that text assumes that $\Xh$ is not forced;
however, the verification of the assumptions given there allows one to
cover that case with the theorems provided in this text.) One uses the
Lyapunov structure $V(u)=U(u)=\LL{u}^2 + \LL{\nabla u}^2$. That case
is also analyzed in \cite{b:Hairer02}. In that reference, the strong
contractive nature is used to get an exponential mixing rate uniform
in the initial data. This is because the time for the initial return
center of the phase space does not depend on the initial state; this
is not the case in the SNS equation. This holds because one can
estimate the time $\tau_C$ uniformly in the initial data. Hence from
\eqref{eq:enterToCReduction}, one sees that the mixing time can be
estimated independent of the initial data. This is made explicitly in
the theorems in \cite{b:Mattingly02}. Another noteworthy feature of
the analysis in \cite{b:Hairer02} is that a change of measure is made
in the low modes to steer all of the modes together only
asymptotically.  In contrast to the presentation given here where the
$\ell$ variable is made to be equal for all moments of time after
$t=1$ and the $h$ variable converges asymptotically. The method in
$\cite{b:Hairer02}$ appears to be simpler to construct while the
method exposed here gives convergence in a slightly stronger topology.

In \cite{b:ELui02} other examples are given, all but the stochastic
Kuramoto-Sivashinsky, fits directly in the framework given here.  The
stochastic Kuramoto-Sivashinsky equation requires localization ideas
not based on a straight forward Lyapunov function.  The details are
explained fully in \cite{b:ELui02}.

\section{True Hypoellipticity and the
  Cascade of Randomness}
\label{sec:hypoelliptic}

It is reasonable to ask if the results given in Theorem
\ref{thm:mainErgodic} or Theorem \ref{thm:abs_ergodic} are sharp. Does
ergodicity require forcing all of the modes below the scale specified
by the balance between energy influx and dissipation ? The assumption
for the second part of Theorem \ref{thm:mainErgodic} is an ellipticity
assumption about the dynamics in the typically unstable directions.
Equivalently viewed from the Memory/Gibbsian dynamics point of view,
it means that the reduced system with memory \eqref{eq:lowModeMemory}
is elliptic.

While there is no complete proof, there are a number of results which
seem to imply that much weaker conditions are sufficient. They all
describe the dynamics in a hypoelliptic setting; the case where all of
the typically unstable degrees of freedom are not forced directly. In
this setting, ergodicity and mixing require that the nonlinearity
transfer the randomness to other degrees of freedom. 

The first result given below proves the ergodicity of an arbitrary Galerkin
approximation of \eqref{eq:vorticity} under very weak assumptions.
Under similar assumptions, the second result says that the full PDE has a
transition density whose finite dimensional marginals have a
density with respect to Lebesgue measure. 

A third result by \cite{b:Romito02Pre} proves the geometric ergodicity
of the Galerkin projections of the three dimensional SNS equations.
This was expected as the structure shares the needed structure with
the two dimensional problem. What was extremely interesting and novel
in that paper was the proof that the system was globally controllable.
A fourth result found in \cite{b:AcrachevSarychev03Pre} shows that the
full two and three dimensional SNS equations are controllable in the
sense that one can steer them so that any finite number of modes take
specified values. This is very similar in spirit to Theorem
\ref{thm:densityPDE} where only projections of the transition measure
are shown to have a density. The techniques used to prove the control
results in \cite{b:Romito02Pre} and \cite{b:AcrachevSarychev03Pre}
seem to use the same important observation. Namely that the off-diagonal
nature of the nonlinearity leaves the system globally consolable even
though its nonlinearity is even powered.  We refer the reader to
\cite{b:Romito02Pre} and \cite{b:AcrachevSarychev03Pre} for the
precise statement of the results.

As we will undertake direct calculations, it is simpler to work in
a real basis of $L^2(\T^2)$. For this reason we switch our forcing
to the form
\begin{align}\label{eq:WNew}
  W(x,t)=\sum_{k\in \mathcal{K}^{\cos}} \sigma_k^{\cos} \cos(k\cdot x)b_k(t) +
  \sum_{k\in \mathcal{K}^{\sin}} \sigma_k^{\sin} \sin(k\cdot x)B_k(t)
\end{align}
where $B_k$ and $b_k$ are independent real Brownian motions with
variance one, $\sigma_k^{\cos}$, $\sigma_k^{\sin}$ are positive real
constants, and $ \mathcal{K}^{\cos}$, $ \mathcal{K}^{\sin}$ are
subsets of $\Z^2_*\eqdef\{ j=(j_1,j_2) \in \Z^2: j_2 \geq 0, |j|>0\}$.
We need only to consider $\Z^2_*$ as the reality of the vorticity
allows one to restrict to wave number in the upper half plane and we
have assumed the absence of a mean flow.  (Note: In
\cite{b:EMattingly00} the sums were restricted too much, however this
does not effect any of the bracket calculations made and the results
hold true.)

We now define two sequences of subsets of $\Z^2$ which capture how
the randomness spreads from one degree of freedom to the next. 
Define $\mathcal{K}_0=\mathcal{Z}_0= \mathcal{K}^{\cos} \cap
\mathcal{K}^{\sin}$. Next define
\begin{multline*}
  \mathcal{Z}_n=   \mathcal{Z}_{n-1} \cup \big\{ k \in \Z^2_* :
  k\in\{\ell+j,\ell-j,j-\ell\}  \mbox{ with } j \in \mathcal{Z}_0, \ell\in
  \mathcal{Z}_{n-1}\\
 \text{ and } \ell^\perp\cdot j \not =0, |j|\not =|\ell| \big\}
\end{multline*}
and fixing some positive integer $N$ define
\begin{multline*}
  \mathcal{K}_n^N=   \mathcal{K}_{n-1}^N \cup \big\{ k \in \Z^2_* :
  k\in\{\ell+j,\ell-j,j-\ell\} \mbox{ with } j, \ell\in
  \mathcal{K}_{n-1}^N\\
 \text{ and } \ell^\perp\cdot j \not =0, |j|\not =|\ell|,
  |\ell-j|<N,|\ell+j| <N  \big\}
\end{multline*}
and finally $\mathcal{Z}_\infty = \cup \mathcal{Z}_n$ and
$\mathcal{K}_\infty^N = \cup \mathcal{K}_n^N$. The two sets track the
cascade of randomness out to the unforced modes. The farther along the chain
which a mode first enters the sequence of sets, the less the random
variation will be felt in that coordinate. 

Theorem \ref{thm:galerkin} below will state its assumptions in terms
of $\mathcal{K}_\infty^N$ whereas Theorem \ref{thm:densityPDE} will
use $\mathcal{Z}_\infty$. It is likely that for a given
$\mathcal{K}^{\cos}$ and $\mathcal{K}^{\sin}$ that
$\mathcal{Z}_\infty=\cup_N\mathcal{K}_\infty^N$ (one direction is
clear) however proof is not immediately obvious. Furthermore, a sketch
of Theorem \ref{thm:densityPDE}, under the same assumptions as Theorem
\ref{thm:galerkin}, is given in \cite{b:MattinglyPardoux03Pre}. Hence,
we do not think there is any real significant difference between the two
sets.

The first result we state gives exponentially mixing for the order $N$
Galerkin approximation of \eqref{eq:vorticity} with forcing of the
form \eqref{eq:WNew} provided an algebraic condition on the wave
numbers forced, given in terms of $\mathcal{K}_\infty^N$, is
satisfied.  By the Galerkin approximation of order $N$, we mean the
finite system of coupled ODEs obtained by setting to zero, for all
time, any Fourier mode with $|k| \geq N$. This approximation returns
us to the setting of standard hypoelliptic SDE in $\RR^d$.  Using a
weak version of H\"oromander's sum of squares theorem (cf.
\cite{b:KusuokaStroock84,b:Norris86,b:Bell95} ), it was shown that the
diffusion has a smooth $C^\infty$ density. Then, using some standard
Markov chain theory for a Harris chain with a Foster-Lyapunov
function, one obtains exponential mixing.
\begin{theorem}\label{thm:galerkin} \cite{b:EMattingly00}
  Consider the  order $N$  Galerkin approximation of the vorticity
  equation \eqref{eq:vorticity}. Assume that $\mathcal{K}_\infty^N=\{k
  \in \Z^2_* : |k|<N\}$. Denoting the solution by $\omega^N$, one
  has the following mixing result.
  
  If $\omega_0^N$ and $\tilde\omega_0^N$ are two initial conditions
  then for any $p \geq 1$ there exist positive constants $B=B(p)$ and
  $\gamma=\gamma(p)$ so that
  \begin{multline*}
   \|P_t(\omega^N_0,\ccdot) - P_t(\tilde\omega^N_0,\ccdot)\|_{TV} 
                \leq \|P_t(\omega^N_0,\ccdot) - P_t(\tilde\omega^N_0,\ccdot)\|_{V_p}
                \\
                \leq B[1+|\omega_0^N|^{2p}+|\tilde\omega_0^N|^{2p}]e^{-\gamma t}
  \end{multline*}
Here $\|\ccdot\|_{TV}$ is the total variation norm on signed measures
and  $\|\ccdot\|_{V(p)}$ is the weighted variational norm defined by 
  \begin{align*}
    \|P_t(\omega^N_0,\ccdot) - P_t(\tilde
    \omega^N_0,\ccdot)\|_{V_p}\eqdef \sup_{\phi\in \mathcal{V}_p}
    \EE\phi(\omega^N(t)) - \EE\phi(\tilde\omega^N(t))
  \end{align*}
with $\mathcal{V}_p=\{ \mbox{measurable $\phi$ with }|\phi(x)|\leq 1
+ |x|^{2p} \}$. Taking $\tilde \omega_0^N$ distributed as the
invariant measure, one obtains exponential convergence to the invariant
measure and uniqueness of the invariant measure.
\end{theorem}

To make this theorem interesting, we need some examples of conditions
on $ \mathcal{K}^{\cos}$ and $\mathcal{K}^{\sin}$ so that it
applies. The following Lemma gives simple conditions under which the
previous and next theorems hold.
\begin{lemma}\cite{b:EMattingly00,b:MattinglyPardoux03Pre}\label{l:fullLattice}
  \begin{itemize}
  \item If $\{ (0,1),(1,1)\} \mbox{ or } \{ (1,0),(1,1)\} \subset
    \mathcal{K}^{\cos} \cap \mathcal{K}^{\sin}$ then
    $\mathcal{Z}_\infty=\Z^2$ and $\mathcal{K}_\infty^N=\{k \in \Z^2_*
    : |k|<N\}$ for any $N$.
  \item    Let $M,K \in \N$ with $M,K > 2$ and $|M-K|> 2$. Then if 
                \begin{gather*}
                        \{ (M+1,0),(M,0),(0,K+1),(0,K) \} \subset \mathcal{Z}_0
                \end{gather*}
    then $\mathcal{Z}_\infty=\Z^2$. If in addition $M,K < N-1$ then
    $\mathcal{K}_\infty^N=\{k \in \Z^2_* : |k|<N\}$.
\end{itemize}
\end{lemma}
This gives only two examples of types of forcing which are
sufficiently distributed to ensure ergodicity. Many others choices are
possible. The author thanks A. Majda and P. Constantin for stimulating
conversations which pushed him to verify the second part of Lemma
\ref{l:fullLattice}. It provides an example of forcing which allows
one to observe both the energy and enstrophy cascade which are present
in two dimensional fluid systems. Of course, the most interesting
question would be to make some qualitative statement connecting this
cascade of probability with the dynamics.

A theorem similar to Theorem \ref{thm:galerkin}, but for the three
dimensional Galerkin approximation, is proven in \cite{b:Romito02Pre}.
There he proves even more; he shows that the system is actually
globally controllable. This very interesting fact hinges on the
observation that because the nonlinearity is off-diagonal in Fourier
space; and hence, the system has the good properties of systems with
odd powered polynomials nonlinearities (see \cite{b:Jurdjevic98}).

Theorem \ref{thm:galerkin} gives a strong indication that a similar
theorem holds for the full PDE; however, a proof currently alludes the
community. The following theorem shows that at least one of the needed
ingredients persists for the full infinite dimensional vorticity equation.

Defining
\begin{align*}
   S_\infty=\mathrm{Span}\left(\big\{ \sin(k\cdot x): k \in \mathcal{Z}_\infty
     \cup \mathcal{Z}^{cos}_0 \big\}
     \cup \big\{ \cos(k\cdot x): k \in \mathcal{Z}_\infty
     \cup \mathcal{Z}^{sin}_0  \big\} \right) 
\end{align*}
we have the following density result for the finite dimensional
marginals of equation \eqref{eq:vorticity}.
\begin{theorem}\label{thm:densityPDE} \cite{b:MattinglyPardoux03Pre}
  For any $t > 0$ and any finite dimensional subspace $S$ of $S_\infty$,
  the law of the orthogonal projection $\Pi_S \omega(t, \cdot)$ of
  $\omega(t, \cdot)$ onto $S$ is absolutely continuous with respect to
  the Lebesgue measure on $S$.
\end{theorem}
This of course is not enough to prove ergodicity. It addresses  only the
first part of Assumption \ref{a:leaves}.

\section{Open Questions}
\label{sec:open}
A number of open questions have been mentioned in the text. Here we
collect them and add a few more.
\begin{enumerate}
\item Extend the ergodic results to the case when all of the
  determining modes are not forced. The results on the ergodicity of
  the Galerkin approximation suggest strongly that full PDE is ergodic
  under weaker assumptions than Theorem \ref{thm:mainErgodic}. Theorem
  \ref{thm:galerkin} gives and indication what the proper assumptions
  should be.  The results on the existence of densities for the
  projection of transition densities and the controllability of a
  finite number of variables gives strong evidence that nothing
  surprising happens in the full PDE.
\item Prove (or disprove) that even when the forcing has spatial
  Fourier modes which decay super--exponentially, the solution
  still decays only exponentially in $|k|$. Prove (or disprove) that
  this decay rate does not fluctuate with time in the stationary state. 
\item Related to the previous: ``What is the natural topology of the
  transition density of the Markov process defined by the SNS ?''
\item Extend Theorem \ref{thm:mainErgodic} to the full space. The case
  of bounded domains in the same as the periodic case. However the
  full space requires some additional ideas, if not completely
  different ones.
\item Understand better the $\nu \rightarrow 0$ limit. In a recent
  preprint \cite{b:Kuksin03Pre} explores this limit for one choice of
  forcing. However, the choice of scaling produces a deterministic
  limit which is the less interesting case and does not correspond to
  the traditional view of turbulence. In all cases, there remain many
  interesting question concerning the structure of the limiting
  solutions and the limit when other types of forcing are used.
\item Make progress in the three dimensional problem. Unless a
  breakthrough is made in the deterministic three dimensional problem,
  this would likely require other methods.  The methods used here
  proceed in a pathwise manner in the high $k$ and, hence, can do no
  better than the deterministic theory. In particular, the estimates
  used to get contraction of the high $k$ are similar to those used to
  prove uniqueness of solutions. Recently Da Prato and Debussche have
  show that by a selection principle one can build a stochastic
  process associated to the 3D problem and that this process under
  certain conditions has a unique invariant measure. Unfortunately the
  conditions on the forcing require it to have algebraic decay in
  $|k|$.
\end{enumerate}

\section{Acknowledgments}
I am indebted to my collaborators Yuri Bakhtin, Weinan E, Ya. Sinai,
Toufic Suidan, Andrew Stuart for both their hand in exploring the
questions described in this note and there advice during its writing.
I also thank Persi Diaconis, Amir Dembo, and George Papanicolaou for
useful discussions when I first worked to understand the remaining
needed estimates during my early days at Stanford.  I also thank
Sandra McBride and BJM for reading sections.  I thank the organizers
of Forges-les-Eaux for the invitation to talk and the opportunity to
publish this expanded version of my lecture. The also author thanks
the Institute for Advanced Study for it hospitality during the year
2002-2003 when the majority of this text was written and the NSF for
its support through Grant DMS-9971087.
\appendix

\section{Comparison of Measures on Path Space}
\label{sec:measurePath} 
Suppose that we have stochastic processes $X^{(i)}(t)$, $i=1,2$ on the
path space $C([0,T],\X)$ where $\X$ is some Hilbert space and $T\in
(0,\infty]$. Furthermore, assume that $X^{(i)}$ satisfies the equation
\begin{equation}
\begin{split}
dX^{(i)}(t)&=f_i(t,X^{(i)}_{[0,t]})dt+gdW(t),\ t \in [0,T]\\
X^{(i)}(0)&=x_0.
\end{split}
\label{2SDDEs}
\end{equation}
Here, for fixed $t$ the functions $f_1$ and $f_2$ map the space
$C_{[0,t]}=C([0,t],\X)$ to $\X$.  By $X_{[0,t]}$ we mean the segment of
the trajectory on $[0,t]$. $W(t)$ is a cylindrical Brownian motion
over a Hilbert space $\Y$ and $g$ is an invertible Hilbert-Schmidt operator
from $\Y \rightarrow \X$.  For any $\mathcal{B}\subset C_{[0,T]}$,
define measures $P^{(i)}_{[0,T]}(\ccdot ; \mathcal{B})$ on the path
space as:
\begin{equation*}
  P^{(i)}_{[0,T]}(A; \mathcal{B}) = 
  P\{X^{(i)}_{[0,T]}\in A\cap\mathcal{B}\},\ \mbox{\rm for} \ A
  \subset C_{[0,T]}. 
\end{equation*}
Define also $D(t,\ccdot)=f_1(t,\ccdot)-f_2(t,\ccdot)$.

In this setting, we have the following result which is a variation on
Lemma B.1 from \cite{b:Mattingly02} and follows quickly from
Girsanov's Theorem. Similar versions of this lemma can be found in
\cite{b:MattinglySuidan03Pre} and \cite{b:BakhtinMattingly03Pre}.
\begin{lemma} \label{l:compareMeasures}  
  Assume there exists a constant $D_*\in(0,\infty)$ such that
  \begin{align}\label{eq:Novikov}
    \exp\left\{\frac12 \int_0^T
      \big|g^{-1}D\big(t,X^{(i)}_{[0,t]}\big)\big|^2_{\Y} dt \right\}
    \one_\mathcal{B}(X^{(i)}_{[0,t]}) < D_*
  \end{align}
  almost surely for $i=1,2$.  Then the measures
  $P^{(1)}_{[0,T]}(\ccdot ; \mathcal{B})$ and $P^{(2)}_{[0,T]}(\ccdot
  ; \mathcal{B})$ are equivalent. In addition for any $p>0$
  \begin{equation*}
    \EE \left[\frac{dP^{(1)}_{[0,T]}(\ccdot ;
      \mathcal{B})}{dP^{(2)}_{[0,T]}(\ccdot ; \mathcal{B})}\right]^p
      \leq D_*^{p(p-1)} \ .
  \end{equation*}
And lastly
\begin{align*}
  \frac12\|P^{(1)}_{[0,T]}(\ccdot ; \mathcal{B})
  -P^{(2)}_{[0,T]}(\ccdot ; \mathcal{B})\|_{TV}\leq  \left( \EE
    \left[\frac{dP^{(1)}_{[0,T]}(\ccdot )}{dP^{(2)}_{[0,T]}(\ccdot )}
        -1 \right]^2 \one_{ \mathcal{B}}\right)^\frac12 \leq \left(
  D_*^2 -1\right)^\frac12
\end{align*}
  \end{lemma}
  \bpf
  Define the auxiliary SDEs
  \begin{align*}
    dY^{(i)}(t) &=f_i\big(t, Y^{(i)}_{[0,t]}\big)
    \one_{\mathcal{B}(t)}(Y^{(i)}_{[0,t]}) dt + gdW(t)
  \end{align*}
  where 
  $
    \mathcal{B}(t)=\{ x \in C_{[0,t]}: \exists \bar{x} \in
    \mathcal{B} \mbox{ such that } x(s)=\bar{x}(s) \mbox{ for $s
    \in[0,t]$}  \} 
  $. 
  Solutions $Y^{(i)}(t)$ to these equation can be constructed as
  \begin{equation*}
    Y^{(i)}(t)=X_i(t)\one_{\{t\le\tau\}}+[gW(t)-gW(\tau)+X^{(i)}(\tau)]\one_{\{t>\tau\}}.   
  \end{equation*}
  Here $\tau=\inf\{s>0: X^{(i)}_{[0,s]}\not\in \mathcal{B}(s)\}$.
  
  Denote
  $D_\mathcal{B}(t,x) =[f_1(t,x) - f_2(t,x)]\one_{\mathcal{B}(t)}(x)$.
  The assumption on $D$ in \eqref{eq:Novikov} and the definition of 
  $\mathcal{B}(t)$ imply that 
  \begin{equation*}
    \exp\left\{\frac12\int_0^T \big| g^{-1}D_\mathcal{B}
    \big(t,X_{[0,t]}\big)\big|^2_{\Y} dt
    \right\}<D_*\quad \mbox{a.s.}
  \end{equation*}
  under both measures $P^{(i)}_{Y[0,t]}$ defining 
  solutions to auxiliary equation with $i=1$
  and $i=2$.
  Hence, Novikov's condition is satisfied for the difference of the drifts
  of the auxiliary equations and 
  the Girsanov's theorem implies that 
  $\frac{dP^{(1)}_{Y[0,t]}}{dP^{(2)}_{Y[0,t]}}(x) = \mathcal{E}(x)$ where
  the Radon--Nikodym derivative evaluated at a trajectory $x$ is defined by
  the stochastic exponent:
  \begin{align*}
    \mathbf{E}(x)=\exp\left\{\int_0^T
     \left\langle g^{-1}D_{\mathcal{B}}(s, x_{[0,s]}), dW(s)\right\rangle_{\Y} - \frac12
      \int_0^T |g^{-1}D_{\mathcal{B}}(s, x_{[0,s]})|^2_{\Y} ds \right\}.
  \end{align*}
 
  Note that restrictions of measures $P^{(i)}_{Y_{[0,t]}}$ on the set
  $\mathcal{B}$ coincide with $P^{(i)}_{[0,t]}(\ccdot;\mathcal{B})$.
  This proves that $P^{(1)}_{[0,t]}(\ccdot,\mathcal{B})$ is absolutely
  continuous with respect to $P^{(2)}_{[0,t]}(\ccdot;\mathcal{B})$.
  The reverse relation follows by symmetry and the proof of
  equivalence is complete.

  To prove the second estimate, notice that 
  \begin{align*}
    (\mathbf{E})^p&=\exp\left\{p\int_0^T
     \left\langle g^{-1}D_{\mathcal{B}}(s, x_{[0,s]}), dW(s)\right\rangle_{\Y} -p \frac12
      \int_0^T |g^{-1}D_{\mathcal{B}}(s, x_{[0,s]})|^2_{\Y} ds \right\}
     \\
     &= \mathbf{E}_p \exp\left( \frac{p^2-p}2  \int_0^T
     |g^{-1}D_{\mathcal{B}}(s, x_{[0,s]})|^2_{\Y} ds  \right) \leq
     \mathbf{E}_p D_*^{p(p-1)}
  \end{align*}
  where $\mathbf{E}_p$ is the martingale defined by
   \begin{align*}
     \mathbf{E}_p=\exp\left\{p\int_0^T \left\langle
         g^{-1}D_{\mathcal{B}}(s, x_{[0,s]}), dW(s)\right\rangle_{\Y}
       - \frac{p^2}2 \int_0^T |g^{-1}D_{\mathcal{B}}(s,
       x_{[0,s]})|^2_{\Y} ds \right\}.
  \end{align*}
  Hence, $\EE \mathbf{E}_p=1$ and in light of the  estimate on
  $\mathbf{E}^p$, the proof is complete. To see the last estimate, use
  the Cauchy-Schwartz inequality to obtain the first inequality. The
  expand the square and use the fact that the Radon-Nikodym derivative
  is a martingale with expectation one to obtain the bound $(
  \EE(\frac{dP^{(1)}}{dP^{(2)}})^2 -1 )^\frac12$. Applying the previous
  estimate to the square gives the result.
  \epf
\section{Coupling Estimates}
\label{sec:couplingEstimatesAppendix}
For any two probability measure $\mu_1$ and $\mu_2$ on a space $\X$, we
can always write them relative to a common measure $\nu$ so that
$d\mu_i=\psi_i d\nu$.  Then we define the measures $(\mu_1 \wedge
\mu_2)(\ccdot)$, $(\mu_1-\mu_2)^+(\ccdot)$, and
$(\mu_2-\mu_1)^+(\ccdot)$ respectively by the densities $(\psi_1
\wedge \psi_2)d\nu$, $(\psi_1 -\psi_2)^+d\nu$, $(\psi_2 -
\psi_1)^+d\nu$ where $a \wedge b=\min(a,b)$ and $(a)^+$ is $a$ if $a$
is positive and zero otherwise.  Notice that $\mu_1= (\mu_1 \wedge
\mu_2) +(\mu_1-\mu_2)^+$. Also observe that if $\|\ccdot\|_{TV}=\int |\psi_1-\psi_2|d\nu$ is the
total variation norm then $\frac12\|\mu_1 - \mu_2\|_{TV}=1-(\mu_1 \wedge
\mu_2)(\X)=(\mu_1-\mu_2)^+(\X)=(\mu_2-\mu_1)^+(\X)$. The proof of the
following lemma can be found in the appendix of \cite{b:Mattingly02}.

\begin{lemma}\label{l:abstractCouplingBound}
  Let $\mu_1$ and $\mu_2$ be two  measures on a space $\X$ with
  $\mu_i(\X) \leq 1$. Assume that $\mu_1$ is
  equivalent to $\mu_2$ and that there exists a constant $C'>0$
  and $p>1$ so that 
  \begin{align*}
    \int_\X\left[\frac{d\mu_1}{d\mu_2}(x)\right]^{p+1} d\mu_2(x)=\int_\X
    \left[\frac{d\mu_1}{d\mu_2}(x)\right]^p d\mu_1(x) < C'
  \end{align*}
  then 
  \begin{align*}
    \int_\X \left| 1\wedge \frac{d\mu_1}{d\mu_2}(x)
    \right| d\mu_2(x) \geq\left[1-\frac1p\right]
    \left(\frac{\mu_1(\X)^p}{pC'}\right)^\frac{1}{p-1}\ .
  \end{align*}
Notice that this lower bound is strictly positive if $\mu_1(\X)>0$ (or
equivalently $\mu_2(\X )>0$).
\end{lemma}

{\tiny
\textbf{Errors/Typos corrected since original version:}
\begin{itemize}
\item[12/03] Fix misplaced ``a'' in Assumption \ref{a:bad2} on p. \pageref{a:bad2}.
Fix missing power of 2 in definitions of
  $\mathcal{E}_\alpha$ and $\mathcal{E}_1$ on p.
  \pageref{er:E1} and  p. \pageref{er:E2}  respectively.
 Correct omitted restriction to $\mathcal{D}(G)$ on p.
  \pageref{er:G1} and p. \pageref{er:G2} and associated rewording of Assumption \ref{a:Lyop}
  on p. \pageref{er:Phi1}.
Clarify assumptions on $\Phi$ on p. \pageref{er:Phi1}. 
Replace cosmetics with cosmetic on p. \pageref{er:cosmetic}. 

\item[2/04] Fix direction of inequality in Lemma \ref{l:returnTimes}.
\end{itemize}
}
\end{document}